\documentclass[12pt]{article}
\usepackage{amsmath,amssymb,latexsym}
\usepackage[a4paper]{geometry}
\usepackage{showkeys}

\newcommand{\eh}{\hfill}
\newlength{\sperr}

\newenvironment{thm}[2]{\begin{sloppypar}{#1 
#2.}\em{}}{\end{sloppypar}}

\newcommand{\proof}{\hspace*{9mm}{\settowidth{\sperr}{\rm 
Proof}\parbox[t]{1.3\sperr}{\rm P\eh r\eh o\eh o\eh f\eh } }}

\newlength{\addro}\newlength{\addrt}
\input{tcilatex}

\begin{document}

\begin{center}
\bigskip {\Large \textbf{On the Lagrangian and Hamiltonian aspects of
infinite -dimensional dynamical systems and their finite-dimensional
reductions}}

\textsc{\ Ya. Prykarpatsky} and \textsc{A.M. Samoilenko}

\textit{*) Institute of Mathematics at the NAS, Kiev 01601, Ukraine, and the
AGH University of Science and Technology, Department of Applied Mathematics,
Krakow 30059 Poland (yarpry\@bnl.gov, sam@imath.kiev.ua)}

15 March 2004
\end{center}

ABSTRACT. {\small A description of \ Lagrangian and Hamiltonian formalisms
naturally arisen from the invariance structure of given nonlinear dynamical
systems on the infinite--dimensional functional manifold is presented. The
basic ideas used to formulate the canonical symplectic structure are
borrowed from the Cartan's theory of differential systems on associated
jet--manifolds. The symmetry structure reduced on the invariant submanifolds
of critical points of some nonlocal Euler--Lagrange functional is described
thoroughly for both differential and differential discrete dynamical
systems. The Hamiltonian representation for a hierarchy of Lax type
equations on a dual space to the Lie algebra of integral-differential
operators with matrix coefficients, extended by evolutions for
eigenfunctions and adjoint eigenfunctions of the corresponding spectral
problems, is obtained via some special Backlund transformation. The
connection of this hierarchy with integrable by Lax spatially
two-dimensional systems is studied.}

\section{Introduction}

\setcounter{equation}{0} \renewcommand{\theequation}{\arabic{section}.%
\arabic{equation}}

One of the fundamental problems in modern theory of infinite-dimensional
dynamical systems is that of an invariant reduction them upon some invariant
submanifolds with enough rich mathematical structures to treat their
properties analytically. The first approaches to these problems were
suggested still at the late times of the preceding century, in the classical
oeuvres by S.Lie, J.Liouville, J.Lagrange, V.R.Hamilton, J.Poisson and
E.Cartan. They introduced at first the important concepts of symmetry,
conservation law, symplectic, Poisson and Hamiltonian structures as well
invariant reduction procedure, which appeared to be extremely useful for
proceeding studies. These notions were widely generalized further by Souriau 
\cite{Sou}, Marsden and Weinstein \cite{Wein,Mars}, Lax \cite{Lax},
Bogoyavlensky and Novikov \cite{Bogo}, as well by many other researchers 
\cite{Pri,Gill,Wah,Kup,Kupe}. It seems worthwhile to mention here also the
recent enough studies in \cite{Adl,Pry,Oew,Oe,Fok,Olv,Mag,Fer}, where the
special reduction methods were proposed for the integrable nonlinear
dynamical systems on both functional and operator manifolds. In the present
paper we describe in detail the reduction procedure for infinite dimensional
dynamical systems upon the invariant set of critical points of some global
invariant functional. The method uses the Cartan's differential-geometric
treating of differential ideals in Grassmann algebra over the associated
jet--manifold. As one of main results, we show also that both the reduced
dynamical systems and their symmetries, generate the Hamiltonian flows on
the invariant critical submanifolds of local and nonlocal functionals with
respect to the canonical symplectic structure upon it. These results are
generalized for the case of differential-difference dynamical systems being
given on discrete infinite-dimensional manifolds. The direct procedure to
construct the invariant Lagrangian functionals for a given apriori Lax-type
integrable dynamical system is presented for both the differential and the
differential-difference cases of the manifold \textit{M.} Some remarks on
the Lagrangian and Hamiltonian formalisms, concerned to infinite-dimensional
dynamical systems with symmetries are given. {\small The Hamiltonian
representation for a hierarchy of Lax type equations on a dual space to the
Lie algebra of integral-differential operators with matrix coefficients,
extended by evolutions for eigenfunctions and adjoint eigenfunctions of the
corresponding spectral problems, is obtained via some special Backlund
transformation. The connection of this hierarchy with integrable by Lax
spatially two-dimensional systems is studied.}

\section{General setting}

\setcounter{equation}{0}

We are interested in treating a given nonlinear dynamical system 
\begin{equation}
du/dt=K[u],
\end{equation}%
with respect to an evolution parameter $t\in \mathbb{R}$ on an
infinite-dimensional functional manifold $M\subset C^{(\infty )}(\mathbb{R};%
\mathbb{R}^{m}),$ possessing two additional ingredients: a homogenous and
autonomous conservation law $\mathcal{L}\in D(M)$ and a number of homogenous
autonomous symmetries $du/dt_{j}=K_{j}[u],\quad j=\overline{1,k},$ with
evolution parameters $t_{j}\in \mathbb{R}$. The dynamical system (2.1) is
not-supposed to be in general Hamiltonian, all the maps $K,K_{j}:M%
\rightarrow T(M),\quad j=\overline{1,k},$ being considered smooth and
well-defined on $M$.

To pose the problem to be discussed further more definitely, let us involve
the jet-manifold $J^{(\infty )}(\mathbb{R};\mathbb{R}^{m})$ locally
isomorphic to the functional manifold $M$ mentioned above. This means the
following: the vector field (2.1) on $M$ is completely equivalent to that on
the jet-manifold $J^{(\infty )}(\mathbb{R};\mathbb{R}^{m})$ via the
representation \cite{Grif,Fil} 
\begin{equation}
(M\ni u\rightarrow K[u])\overset{jet}{\longrightarrow (}(K(u,u^{(1)},\dots
,u^{(n+1)})\leftarrow (x;\,u,u^{(1)},\dots ,u^{(\infty )})\in J^{(\infty )}(%
\mathbb{R};\mathbb{R}^{m})),
\end{equation}%
where $n\in \mathbf{Z}_{+}$ is fixed, $x\in \mathbb{R}$ is the function
parameter of the jet-bundle $J^{(\infty )}(\mathbb{R};\mathbb{R}^{m})\overset%
{\pi }{\longrightarrow }\mathbb{R},$ and $\pi $ is the usual projection on
the base $\mathbb{R}$. Let us allow also that the smooth functional $%
\mathcal{L}\in D(M)$ is a conservation law of the dynamical system (2.1),
that is $d\mathcal{L}/dt=0$ along orbits of (2.1) for all $t\in \mathbb{R}.$
Due to the jet-representation (2.2) we can write the density of the
functional $\mathcal{L}\in D(M)$ in the following form: 
\begin{equation}
\mathcal{L}=\int_{\mathbb{R}}dx\mathcal{L}[u],
\end{equation}%
with $\mathbb{R}\mathbf{\times }\mathbb{R}^{m}\ni \lbrack x;u]\overset{jet}{%
\longrightarrow }(x;u,u^{(1)},\dots ,u^{(N+1)})\in J^{(N+1)}(\mathbb{R};%
\mathbb{R}^{m})$ being the standard jet-mapping and a number $N\in \mathbf{Z}%
_{+}$ fixed. Besides, the functional (2.3) will be assumed to be
non-degenerate in the sense that Hessian of $\mathcal{L}:J^{(N+1)}(\mathbb{R}%
;\mathbb{R}^{m})\!\!\rightarrow \mathbb{R}$ has nonvanishing determinant: $%
det\Vert \frac{\partial ^{2}\mathcal{L}(u,u^{(1)},\dots ,u^{(N+1)})}{%
\partial u^{(N+1)}\partial u^{(N+1)}}\Vert \neq 0.$

\section{Lagrangian reduction}

\setcounter{equation}{0}

Consider now the set of critical points $M_{n}\subset M$ of the functional $%
\mathcal{L}\in D(M)$: 
\begin{equation}
M_{N}=\{u\in M:\mathop{\rm grad}\,\mathcal{L}[u]=0\}
\end{equation}%
where, due to (2.2), $\mathop{\rm grad}\,\mathcal{L}[u]:=\delta \mathcal{L}%
(u,\dots ,u^{(N+1)})/\delta u$ -- the Euler variational derivative. As
proved by Lax \cite{Lax}, the manifold $M_{N}\subset M$ is smoothly imbedded
well-defined one due to the condition $Hess\mathcal{L}\neq 0$. Besides, the
manifold $M_{N}$ appears to be invariant with respect to the dynamical
system (2.1). This means in particular that the Lie-derivative of any vector
field $X:M\rightarrow T(M),$ tangent to the manifold $M_{N}$ with respect to
the vector field (2.1), is again tangent to $M_{N},$ that isthe followung
implication%
\begin{equation}
X[u]\in T_{u}(M_{N})\Rightarrow \lbrack K,X][u]\in T_{u}(M_{N})
\end{equation}%
holds for all $u\in M_{N}.$ Here we are in a position to begin with a study
of the intrinsic structure of the manifold $M_{N}\subset M$ within the
geometric Cartan's theory developed on the jet-manifold $J^{(\infty )}(%
\mathbb{R};\mathbb{R}^{m})$ \cite{Fil,Lax,Bry,St}. Let us define an ideal $%
I(\xi )\subset \Lambda (J^{(\infty )})$, generated by the vector one-forms $%
\xi _{j}^{(1)}=du^{(j)}-u^{(j+1)}dx,\quad j\in \mathbf{Z}_{+},$ which are
canceled the vector field $d/dx$ on the jet-manifold $J^{(\infty )}(\mathbb{R%
};\mathbb{R}^{m}):$ 
\begin{equation}
i_{\frac{d}{dx}}\xi _{j}^{(1)}=0,\qquad j\in \mathbf{Z}_{+},
\end{equation}%
where $x\in \mathbb{R}$ belongs to the jet-bundle base$,\quad i_{\frac{d}{dx}%
}$ -- the intrinsic derivative along the vector field 
\begin{equation*}
\frac{d}{dx}=\frac{\partial }{\partial x}+\sum_{j\in \mathbf{Z}%
_{+}}<u^{(j+1)},\frac{\partial }{\partial u^{(j)}}>,
\end{equation*}%
where $\langle .,.\rangle $ is the standard scalar product in $\mathbb{R}%
^{m} $. The vector field (2.1) on the jet-manifold $J^{(\infty )}(\mathbb{R};%
\mathbb{R})$ has the analogous representation: 
\begin{equation}
\frac{d}{dt}=\frac{\partial }{\partial t}+\sum_{j=\mathbf{Z}_{+}}<K^{(j)},%
\frac{\partial }{\partial u^{(j)}}>,
\end{equation}%
where, by definition, $K^{(j)}:=\frac{d^{j}}{dx^{j}}K,\,j=\mathbf{Z}_{+}.$
The problem arises: how to build the intrinsic variables on the manifold $%
M_{N}\subset M$ from the jet-manifold coordinates on $J^{(\infty )}(\mathbb{R%
};\mathbb{R}^{m}))$? \ 

To proceed to the solution of the problem above, let us study the 1-form $d%
\mathcal{L}=\Lambda ^{1}(J^{(\infty )}(\mathbb{R};\mathbb{R}))$ as one
defined on the submanifold $M_{N}\subset M.$ We have the following chain of
identities in the Grassmann subalgebra $\Lambda (J^{(2N+2)}(\mathbb{R};%
\mathbb{R}^{m})):$ 
\begin{eqnarray}
d\mathcal{L} &=&d(i_{\frac{d}{dx}}\mathcal{L}dx)=di_{\frac{d}{dx}}(\mathcal{L%
}dx+\sum_{j=0}^{N}\langle p_{j},\mathbb{R}\rangle ) \\
&=&(di_{\frac{d}{dx}}+i_{\frac{d}{dx}}d)\cdot (\mathcal{L}%
dx+\sum_{j=0}^{N}\langle p_{j},\xi _{j}^{(1)}\rangle )-i_{\frac{d}{dx}}d(%
\mathcal{L}dx+\sum_{j=0}^{N}\langle p_{j},\xi _{j}^{(1)}\rangle ),  \notag
\end{eqnarray}%
where $p_{j}:J^{(2N+2)}(\mathbb{R};\mathbb{R}^{m})\rightarrow \mathbb{R}%
^{m},j=\overline{0,N},$ are some until not definite vector-functions.
Requiring now that 2-form $d(\mathcal{L}dx+\sum_{j=0}^{N}\langle p_{j},\xi
_{j}^{(1)}\rangle )$ 
do not depend on differentials $du^{(j)},j=\overline{1,N+1},$ 
that is%
\begin{equation}
i_{\frac{\partial }{\partial u^{(j)}}}(d\mathcal{L}\wedge
dx+\sum_{k=0}^{N}\langle dp_{k}\wedge \xi _{j}^{(1)}\rangle )=0,
\end{equation}%
one can determine thus the vector-functions $p_{j}=\mathbb{R}^{m},\quad j=%
\overline{0,N}.$ As a result we obtain the following simple recurrent
relations: 
\begin{equation}
\frac{dp_{j}}{dx}+p_{j-1}=\frac{\partial \mathcal{L}}{\partial u^{(j)}}
\end{equation}%
for $j=\overline{1,N+1}$, setting $p_{-1}=0=p_{N+1}$ by definition. The
unique solution to (3.7) is made by the following expressions, $j=\overline{%
0,N}:$ 
\begin{equation}
p_{j}=\sum_{k=0}^{N}{(-1)}^{k}\frac{d^{k}}{dx^{k}}\frac{\partial \mathcal{L}%
}{\partial u^{(j+k+1)}}.
\end{equation}%
Thereby we have got, owing (3.5) and (3.6), the following final
representation for the differential $d\mathcal{L}$: 
\begin{eqnarray}
d\mathcal{L} &=&\frac{d}{dx}[\mathcal{L}-\sum_{j=0}^{N}\langle
p_{j},u^{(j+1)}\rangle ]dx-\langle \mathop{\rm grad}\,\mathcal{L}%
[u],u^{(1)}\rangle dx+  \notag \\
&&+\frac{d}{dx}\biggl(\sum_{j=0}^{N}\left\langle p_{j},du^{(j)}\right\rangle
]\biggr)+\langle \mathop{\rm grad}\,\mathcal{L}[u],du\rangle ,
\end{eqnarray}%
with $\frac{d}{dx}:=di_{\frac{d}{dx}}+i_{\frac{d}{dx}}d$ being the
Lie-derivative along the vector field $\frac{d}{dx},$ and\linebreak $%
\mathop{\rm grad}\,\mathcal{L}[u]:=\delta \mathcal{L}/\delta u$ as it was
mentioned above in the chapter 2. Below we intend to treat the
representation (3.9) on the topic of a symplectic structure arisen from the
above analysis on the invariant submanifold $M_{N}\subset M$.

\section{Symplectic analysis and Hamiltonian formulation}

\setcounter{equation}{0}

Let us put into the expression (3.9) the condition $\mathop{\rm grad}\,%
\mathcal{L}[u]=0$ for all $u=M_{N}.$ Then the following equality is
satisfied: 
\begin{equation}
d\mathcal{L}=\frac{d}{dx}\alpha ^{(1)},\qquad \alpha
^{(1)}=\sum_{j=0}^{N}\langle p_{j},du^{(j)}\rangle ,
\end{equation}%
since the function $h^{(x)}:=\sum_{j=0}^{N}\langle p_{j},u^{(j+1)}\rangle -%
\mathcal{L}(u,\dots ,u^{(N+1)})$ satisfies the condition $%
dh^{(x)}/dx=-\langle \mathop{\rm grad}\,\mathcal{L}[u],u^{(1)}\rangle $ for
all $x=\mathbb{R},$ owing to the relations (3.7). Taking now the external
derivative of (4.1), we obtain that 
\begin{equation}
\frac{d}{dx}\Omega ^{(2)}=0,\qquad \Omega ^{(2)}=d\alpha ^{(1)},
\end{equation}%
where we have used the well known identity $d\cdot \frac{d}{dx}=\frac{d}{dx}%
\cdot d$. From (4.2) we can conclude that the vector field $d/dx$ on the
submanifold $M_{N}\subset M$ is Hamiltonian with respect to the canonical
symplectic structure $\Omega ^{(2)}=\sum_{j=0}^{N}\langle dp_{j}\wedge
du^{(j)}\rangle .$ It is a very simple exercise to state that the function $%
h^{(x)}:J^{(2N+2)}(\mathbb{R};\mathbb{R}^{m})\rightarrow \mathbb{R}$ defined
above is playing a role of the corresponding Hamiltonian function for the
vector field $d/dx$ on $M_{N},$ i.e. the equation 
\begin{equation}
dh^{(x)}=-i_{\frac{d}{dx}}{\Omega }^{(2)}
\end{equation}%
holds on $M_{N}.$ Therefore, we have got the following theorem.

$\QTR{sl}{Theorem}$ $\QTR{sl}{1.}$\textit{\ The critical submanifold M}$%
_{N}\subset $\textit{M defined by (3.1) \ for a given non degenerate smooth
functional L}$=$\textit{D(M)}$\subset $\textit{D(J}$^{(N+1)}$\textit{(}%
\textbf{R;R}$^{m}$\textit{)), being imbedded into the jet-manifold J}$%
^{(\infty )}$\textit{(}\textbf{R;R}$^{m}$\textit{), carries the\ canonical
symplectic structure subject to which the induced vector field \ d/dx on M}$%
_{N}$\textit{\ is Hamiltonian.}

The theorem analogous to the that above was stated before via different
manners by many authors \cite{Pri,Gelf}. Our derivation presented here is as
much simpler as constructive, giving rise to all ingredients of symplectic
theory, stemming from imbedding the invariant submanifold $M_{N} $ into the
jet-manifold.

Now we are going to proceed further to studying the vector field (2.1) on
the manifold $M_{N}\subset M$ endowed with the symplectic structure $\Omega
^{(2)}=\Lambda ^{2}(J^{(N+1)}(\mathbb{R};\mathbb{R}^{m})),$ having built via
the formula (4.2).

We have the following implicating identities: 
\begin{eqnarray}
\frac{d\mathcal{L}}{dt} &=&0\quad \Rightarrow \quad \langle \mathop{\rm grad}%
\,\mathcal{L}[u],K[u]\rangle =-\frac{dh^{(t)}}{dx},  \notag \\
\frac{d\mathcal{L}}{dx} &=&0\quad \Rightarrow \quad \langle \mathop{\rm grad}%
\,\mathcal{L}[u],\frac{du}{dx}\rangle =-\frac{dh^{(x)}}{dx},
\end{eqnarray}%
where functions $h^{(t)}$ and $h^{(x)}$ serve as corresponding Hamiltonian
ones for the vector fields $d/dt$ and $d/dx$. This means in part that the
following equations hold: 
\begin{equation}
dh^{(x)}=-i_{\frac{d}{dx}}\Omega ^{(2)},\quad dh^{(t)}=-i_{\frac{d}{dt}%
}\Omega ^{(2)}.
\end{equation}%
To prove the above statement (4.5), we shall build the following quantities
(for the vector field $d/dx$ at first): 
\begin{equation}
d\cdot i_{\frac{d}{dx}}\langle \mathop{\rm grad}\,\mathcal{L}[u],du\rangle =-%
\frac{d}{dx}(dh^{(x)})
\end{equation}%
stemming from (4.4), and 
\begin{equation}
i_{\frac{d}{dx}}\cdot d\langle \mathop{\rm grad}\,\mathcal{L}[u],du\rangle =-%
\frac{d}{dx}\left( i_{\frac{d}{dx}}\Omega ^{(2)}\right)
\end{equation}%
stemming from (3.9), where we used the preliminary used the evident identity 
$[i_{\frac{d}{dx}},\frac{d}{dx}]=0.$ Adding now the expression (4.6) and
(4.7) entails the following one: 
\begin{equation}
\frac{d}{dx}\langle \mathop{\rm grad}\,\mathcal{L}[u],du\rangle =-\frac{d}{dx%
}(dh^{(x)}+i_{\frac{d}{dx}}\Omega ^{(2)})
\end{equation}%
for all $x=\mathbb{R}$ and $u=M$. Since $\mathop{\rm grad}\,\mathcal{L}[u]=0$
for all $u=M_{N}$, we obtain from (4.8) that the first equality in (4.5) is
valid in the case of the vector field $d/dx$ reduced on $M_{N}$. The
analogous procedure fits also for the vector field $d/dt$ reduced on the
manifold $M_{N}\subset M$. The even difference of the procedure above stems
from the condition on vector fields $d/dt$ and $d/dx$ to be commutative, $%
[d/dt,d/dx]=0$ what engenders the needed identity $[i_{\frac{d}{dt}},\frac{d%
}{dx}]\equiv i_{[\frac{d}{dt},\frac{d}{dx}]}=0$ as a simple consequence of
the procedure considered above. There upon we have stated the validity of
equations (4.5) completely.$\rhd $

\textbf{Theorem 2. }Dynamical systems $d/dt$ and $d/dx$ reduced on the
invariant submanifold $M_{N}\subset M$ (3.1) are Hamiltonian ones with the
corresponding Hamiltonian functions built from the equations (4.4) in the
unique way.

By the way, we have stated also that the Hamiltonian functions $h^{(x)}$ and 
$h^{(t)}$ on the submanifold $M_{N}\subset M$ are commuting with each other,
that is $\{h^{(t)},h^{(x)}\}=0$ where $\{\cdot ,\cdot \}$ is the Poissin
structure on $D(M_{N})$ \ corresponding to the symplectic structure (4.2).
This indeed follows from the equalities (4.4), since $\{h^{(t)},h^{(x)}\}=%
\frac{dh^{(x)}}{dt}=-\frac{dh^{(t)}}{dx}\equiv 0$ upon the manifold $%
M_{N}\subset M.$

\section{Symmetry invariance}

\setcounter{equation}{0}

Let us consider now any vector field $K_{j}:M\rightarrow T(M),\quad j=%
\overline{1,k,}$ being symmetry fields related with the given vector field
(2.1), i.e. $[K,K_{j}]=0,$ $=\overline{1,k}.$ As the conservation law $%
\mathcal{L}\in D(M)$ for the vector field\ (2.1) has not to be that for the
vector fields\ $K_{j},\ j=\overline{1,k},$ the submanifold $M_{N}\subset M$
has not to be invariant also with respect to these vector fields. Therefore,
if a vector field\ $X\in T(M_{N}),$ the vector field\ $[K_{j},X]\notin
T(M_{N})$ in general, if $\frac{d}{dt_{j}},\ j=\overline{1,k},$ are chosen
as symmetries of (2.1). Let us consider the following identity for some
there existing function $\tilde{h}_{j}:J^{(2N+2)}(\mathbb{R};\mathbb{R}%
^{m})\rightarrow \mathbb{R},\ j=\overline{1,k},$ stemming from the
conditions $[\frac{d}{dt},\frac{d}{dt_{j}}]=0,\ j=\overline{1,k},$ on $M:$ 
\begin{equation}
\frac{d}{dt}i_{\frac{d}{dt_{j}}}<\mathop{\rm grad}\,\mathcal{L}[u],du>=-%
\frac{d\tilde{h}_{j}[u]}{dx}.
\end{equation}

\textsl{Lemma 1}$\mathbf{.}$\textit{The functions }$\tilde{h}_{j}[u],$%
\textit{\ }$\ j=\overline{1,k},$\textit{\ reduced on the invariant
submanifold }$M_{N}\subset M$\textit{\ turn into constant ones. These
constants can be chosen obviously as zero ones.}

\textit{Proof.} We have: $[\frac{d}{dt},i_{\frac{d}{dt_{j}}}]=0,$ $j=%
\overline{1,k},$ there upon 
\begin{equation*}
i_{\frac{d}{dt_{j}}}(i_{\frac{d}{dt}}d+di_{\frac{d}{dt}})<\mathop{\rm grad}\,%
\mathcal{L}[u],du>=-\frac{d\tilde{h}_{j}}{dx}\quad \Rightarrow
\end{equation*}%
\begin{equation*}
\Rightarrow i_{\frac{d}{dt_{j}}}\biggl(i_{\frac{d}{dt}}d<\mathop{\rm grad}\,%
\mathcal{L}[u],du>+di_{\frac{d}{dt}}<\mathop{\rm grad}\,\mathcal{L}[u],du>%
\biggr)
\end{equation*}%
\begin{equation*}
=-i_{\frac{d}{dt_{j}}}i_{\frac{d}{dt}}\frac{d}{dx}\Omega ^{(2)}-i_{\frac{d}{%
dt_{j}}}\frac{d}{dx}(dh^{(+)})=-\frac{d}{dx}i_{\frac{d}{dt_{j}}}(i_{\frac{d}{%
dt}}\Omega ^{(2)}+dh^{(t)})=-\frac{d\tilde{h}_{j}}{dx}.
\end{equation*}

Whence we obtain that upon the whole jet-manifold $M\subset J^{\infty }(%
\mathbb{R};\mathbb{R}^{m})$ the following identities hold: 
\begin{equation}
i_{\frac{d}{dt_{j}}}(dh^{(t)}+i_{\frac{d}{dt}}\Omega ^{(2)})=\tilde{h}_{j}%
\text{ }
\end{equation}%
for all $j=\overline{1,k}.$ Since upon the submanifold $M_{N}\subset M\quad $%
onr has $i_{\frac{d}{dt}}\Omega ^{(2)}=-dh^{(t)},$ we find that $\tilde{h}%
_{j}\equiv 0,$ $\ j=\overline{1,k},$ that proves the lemma.$\rhd $

\textit{Note 1. }The result above could be stated also using the standard
functional operator calculus of [8]. Indeed, 
\begin{eqnarray}
\lefteqn{\frac{d}{dt}(i_{\frac{d}{dt_{j}}}<grad\mathcal{L}[u],du>)} \\
&=&\frac{d}{dt}<\mathop{\rm grad}{\cal L}[u],K_{j}[u]>  \notag \\
&=&\langle \frac{d}{dt}\mathop{\rm grad}\,\mathcal{L}[u],K_{j}[u]\rangle +<%
\mathop{\rm grad}\,\mathcal{L}[u],\frac{d}{dt}K_{j}[u]>  \notag \\
&=&\langle -{K^{\prime }}^{\ast }\mathop{\rm grad}\,\mathcal{L}%
[u],K_{j}[u]\rangle +<\mathop{\rm grad}\,\mathcal{L}[u],K_{j}^{\prime }\cdot
K[u]\rangle  \notag \\
&=&-\langle {K^{\prime }}^{\ast }\mathop{\rm grad}\,\mathcal{L}%
[u],K_{j}[u]\rangle +<\mathop{\rm grad}\,\mathcal{L}[u],K^{\prime }\cdot
K_{j}[u]\rangle  \notag \\
&=&-\frac{d}{dx}\mathcal{H}[\mathop{\rm grad}\,\mathcal{L}[u],K_{j}[u]]=-%
\frac{d\tilde{h}_{j}[u]}{dx}.  \notag
\end{eqnarray}%
Here a bilinear form $\mathcal{H}\{\cdot .\cdot \}$ is found via the usual
definition of the adjoined operator ${K^{\prime }}^{\ast }$ for a given
operator $K^{\prime }:L_{2}\rightarrow L_{2}$ with respect to the natural
scalar bracket $(\cdot ,\cdot )$ 
\begin{equation}
({K^{\prime }}^{\ast }a,b):=(a,K{^{\prime }}b),\qquad (a,b):=\int_{\mathbb{R}%
}dx\langle a,b\rangle ,
\end{equation}%
whence we simply obtain: 
\begin{equation}
\langle K{^{\prime }}^{\ast }a,b\rangle -\langle a,K^{\prime }b\rangle =d%
\mathcal{H}[a,b]/dx
\end{equation}%
for all $a,b\in L_{2}.$ Therefore, we can identify $\tilde{h}_{j}[u]=%
\mathcal{H}[\mathop{\rm grad}\,\mathcal{L}[u],K_{j}[u]]$ for all $u\in M,$ $%
j=\overline{1,k}.$ If $u\in M_{N},$ we therewith obtain that $\tilde{h}%
_{j}[u]\equiv 0,\quad j=\overline{1,k},$ that was needed to prove.$\rhd $

As a result of the Lemma proven above one gets the following: the functions $%
\tilde{h}_{j}[u],j=\overline{1,k},$ can not serve as nontrivial Hamiltonian
ones for the dynamical systems $d/dt_{j},j=\overline{1,k},$ upon the
submanifold $M_{N}\subset M.$ To overcome this difficulty we assume the
invariant submanifold $M_{N}\subset M$ to possess some additional symmetries 
$d/dt_{j},j=\overline{1,k},$ which satisfy the following characteristic
criterion: $L_{\frac{d}{dt_{j}}}\mathop{\rm grad}\,\mathcal{L}[u]=0,j=%
\overline{1,k},$ for all $u\in M_{N}.$ This means that for $j=\overline{1,k}$
\begin{equation}
L_{\frac{d}{dt_{j}}}\mathop{\rm grad}\,\mathcal{L}[u]=G_{j}(\mathop{\rm grad}%
\,\mathcal{L}[u]),
\end{equation}%
where $G_{j}(\cdot ),j=\overline{1,k},$ are some linear vector-valued \
functionals on $T^{\ast }(M).$ Otherwise, equations (5.6) are equivalent to
the following: 
\begin{equation}
i_{\frac{d}{dt_{j}}}<\mathop{\rm grad}\,\mathcal{L}[u],du\rangle
=-dh_{j}[u]/dx+g_{j}(\mathop{\rm grad}\,\mathcal{L}[u]),
\end{equation}%
where $g_{j}(\cdot ),j=\overline{1,k},$are some scalar linear functionals on 
$T^{\ast }(M).$ From (3.9) and (5.7) we therewith find that for all $j=%
\overline{1,k}$ 
\begin{equation}
L_{\frac{d}{dt_{j}}}<\mathop{\rm grad}\,\mathcal{L}[u],du\rangle =-\frac{d}{%
dx}(dh_{j}[u]+i_{\frac{d}{dt_{j}}}\Omega ^{(2)})+dg_{j}(\mathop{\rm grad}\,%
\mathcal{L}[u]).
\end{equation}%
If we put now $u\in M_{N},$ that is $\mathop{\rm grad}\,\mathcal{L}[u]=0,$
we immediately will find the following: for all $j=\overline{1,k},$ 
\begin{equation}
dh_{j}[u]+i_{\frac{d}{dt_{j}}}\Omega ^{(2)}=0.
\end{equation}%
Whence we can make a conclusion that the vector fields\ $d/dt_{j},j=%
\overline{1,k},$ are Hamiltonian too on the submanifold $M_{N}\subset M.$
Since $dh_{j}/dx=\{h^{(x)},h_{j}\}=0,j=\overline{1,k},$ on the manifold $%
M_{N},$ we therewith obtain that $dh^{(x)}/dt_{j}=0,j=\overline{1,k}.$ This
is also an obvious corollary of the commutativity $[d/dt_{j},d/dx]=0,$ $\ j=%
\overline{1,k},$ for all $x,t_{j}\in \mathbb{R}$ on the whole manifold $M.$
Indeed, in general case we have the identity $\{h^{(x)},h_{j}\}=i_{[\frac{d}{%
dt_{j}},\frac{d}{dx}]}\Omega ^{(2)},$ whence the equalities $%
\{h^{(x)},h_{j}\}\equiv 0$ hold for all$\ \ j=\overline{1,k}$ \ on the
submanifold $M_{N}\subset M,$ since $[\frac{d}{dt_{j}},\frac{d}{dx}]=0$ on $%
M_{N}$ due to (5.8). The analysis fulfilled above makes it possible to treat
given vector fields\ $d/dt_{j},j=\overline{1,k},$ satisfying either
conditions (5.6) or conditions (5.7) on the canonical symplectic
jet-submanifold $M_{N}\subset M$ analytically as Hamiltonian systems.

\section{Liouville integrability}

\setcounter{equation}{0}

Now we suppose that the vector field\ $d/dt_{j},j=\overline{1,k},$ are
independent and commutative both to each other on the jet-submanifold $%
M_{N}\subset M$ and with vector fields $d/dt$ and $d/dx$ on the manifold $M.$
Besides, the submanifold $M_{N}\subset M$ is assumed to be compact and
smoothly imbedded into the jet-manifold $J^{(\infty )}(\mathbb{R};\mathbb{R}%
^{m}).$ If the dimension $dimM_{N}=2k+4,$ due to the Liouville theorem \cite%
{Pri,Gelf} the dynamical systems $d/dx$ and $d/dt$ are Hamiltonian and
integrable by quadratures on the submanifold $M_{N}\subset M.$ This is the
case for all Lax-integrable nonlinear dynamical systems of the Korteweg-de
Vries type \cite{Lax,Bogo,Gelf,Pri} on spatially one-dimensonal functional
manifolds.

\section{Discrete dynamical systems}

Let us be given a differential discrete smooth dynamical system 
\begin{equation}
du_{n}/dt=K_{n}[u]
\end{equation}%
with respect to a continuous evolution parameter $t\in \mathbb{R}$ on the
infinite-dimensional discrete manifold $M\subset L_{2}(\mathbf{Z};\mathbb{R}%
^{m})$ infinite vector-sequences under the condition of rapid decrease in $%
n\in \mathbf{Z}:\sup_{n\in \mathbf{Z}}{|n|}^{k}\Vert u_{n}\Vert _{\mathbb{R}%
^{m}}<\infty $ for all $k\in \mathbf{Z}_{+}$ at each point $u=(\dots
,u_{n},u_{n+1},\dots )\in M,$ where $u_{n}\in \mathbb{R}^{m},\quad n\in 
\mathbf{Z}.$

Assume further that the dynamical system (7.1) possesses a conservation law $%
\mathcal{L}\in D(M),$ that is $d\mathcal{L}/dt=0$ along orbits of (7.1). Via
the standard operatorial analysis \ one gets from (3.5) the variational
derivative of a functional $\mathcal{L}:=\sum_{n\in \mathbf{Z}}\mathcal{L}%
_{n}[u]:$ 
\begin{equation}
\mathop{\rm grad}\,\mathcal{L}_{n}:=\delta \mathcal{L}[u]/\delta u_{n}={%
\mathcal{L}_{n}^{\prime }}^{\ast }[u]\cdot 1,
\end{equation}%
where the last right-hand operation of multiplying by unity is to be
fulfilled by component.

\textsl{Lemma 2}\textbf{.}\textit{\ Let }$\Lambda (M)$\textit{\ be the
infinite-dimensional Grassmannian algebra on the manifold }$M$\textit{; then
the differential }$dL_{n}[u]\in \Lambda ^{1}(M)$\textit{\ enjoys the
following reduced representation: } 
\begin{equation}
d\mathcal{L}_{n}[u]=<\mathop{\rm grad}\,\mathcal{L}_{n},du_{n}\rangle
+d/dn\,\alpha _{n}^{(1)}[u],
\end{equation}%
\textit{where a one-form }$\alpha _{n}^{(1)}[u]\in \Lambda ^{1}(M)$\textit{\
is determined in unique way, }$\langle \cdot ,\cdot \rangle $ \textit{is the
usual scalar product in }$R^{m}$\textit{\ and }$d/dn=\Delta -1,\quad \Delta $
\textit{is the usual shift operator}.

\textit{Proof.} By definition we obtain for the external differential $d%
\mathcal{L}_{n}[u]$ the following chain of representations for each $n\in 
\mathbf{Z}:$ 
\begin{equation}
d\mathcal{L}_{n}[u]=\sum_{k=0}^{N}<\frac{\partial \mathcal{L}_{n}[u]}{%
\partial u_{n+k}},du_{n+k}>  \notag
\end{equation}%
\begin{equation*}
=\sum_{k=0}^{N}\sum_{s=0}^{k}\frac{d}{dn}<\frac{\partial \mathcal{L}_{n-s}[u]%
}{\partial u_{n+k-s}},du_{n+k-s}>+\sum_{k=0}^{N}<\frac{\partial \mathcal{L}%
_{n-k}[u]}{\partial u_{n}},du_{n}>
\end{equation*}%
\begin{equation*}
=\frac{d}{dn}\sum_{k=0}^{N}\sum_{s=0}^{k}<\frac{\partial \mathcal{L}_{n-s}[u]%
}{\partial u_{n+k-s}},du_{n+k-s}>+\sum_{k=-N}^{N}\mathcal{4}^{-k}<\frac{%
\partial \mathcal{L}_{n}[u]}{\partial u_{n+k}},du_{n}>
\end{equation*}%
\begin{equation*}
=\frac{d}{dn}\sum_{k=0}^{N}\sum_{s=0}^{k}<\frac{\partial \mathcal{L}_{n-s}[u]%
}{\partial u_{n+k-s}},du_{n+k-s}>+<\mathcal{L}_{n}^{^{,\prime }\ast
}.1,du_{n}>
\end{equation*}%
\begin{equation}
=\frac{d}{dn}\alpha _{n}^{(1)}[u]+<\mathop{\rm grad}\,\mathcal{L}%
_{n}[u],du_{n}\rangle ,
\end{equation}%
where $N\in \mathbf{Z}_{+}$ is the fixed number depending on the jet-form of
a functional $\mathcal{L}\in D(M),$ 
\begin{equation*}
\alpha _{n}^{(1)}[u]=\sum_{k=0}^{N}\sum_{s=0}^{k}<\frac{\partial \mathcal{L}%
_{n-s}[u]}{\partial u_{n+k-s}},du_{n+k-s}>
\end{equation*}%
%
%
%
%
%
%
%
%
%
%
%
%
%
%
%
%
%
%
%
%
%
%
%
\begin{equation}
=\sum_{k=0}^{N}\sum_{j=0}^{k}<\frac{\partial \mathcal{L}_{n+j-k}[u]}{%
\partial u_{n+j}},du_{n+j}>,
\end{equation}%
and 
\begin{equation*}
\mathop{\rm grad}\,\mathcal{L}_{n}[u]={\mathcal{L}_{n}^{\prime }}^{\ast
}\cdot 1=\sum_{k=0}^{N}\frac{\partial \mathcal{L}_{n-k}[u]}{\partial u_{n}}.
\end{equation*}%
The latter equality in (7.4) proves the lemma 2 completely.$\rhd $

The above proved representation (7.3) gives rise to the following stationary
problem being posed on the manifold $M:$ 
\begin{equation}
M_N=\{u\in M:\mathop{\rm grad}\,\mathcal{L}_n=0\}
\end{equation}
for all $n\in \mathbf{Z},$ where by definition $det\left\|\frac{\partial^2%
\mathcal{L}n[u]} {\partial u_{N+1}\partial u_{N+1}}\right\|=0.$ In virtue of
(7.3) we obtain the validity of the following theorem.

\textsl{Theorem}\textbf{\ 3.} \textit{The finite-dimensional Lagrangian
submanifold }$M_{N}\subset M$\textit{\ defined by (7.6), is a symplectic one
with the canonical symplectic structure }$\Omega _{n}^{(2)}=d\alpha
_{n}^{(1)}$\textit{\ being independent of the discrete variable }$n\in 
\mathbf{Z}.$

\textit{Proof.} From (7.3) we have that on the manifold $M_{N}\subset M\quad
d\mathcal{L}_{n}[u]=d/dn(\alpha _{n}^{(1)}[u]),$ whence for all $n\in 
\mathbf{Z}\quad d/dn(\Omega _{n}^{(2)})=0.$ This means obviously, that $%
\Omega _{n+1}^{(2)}=\Omega _{n}^{(2)}$ for all $n\in \mathbf{Z},$ or
equivalently, the 2-form $\Omega _{n}^{(2)}$ is not depending on the
discrete variable $n\in \mathbf{Z}.$ As the 2-form $\Omega
_{n}^{(2)}:=d\alpha _{n}^{(1)}$ by definition, this form is chosen to be a
symplectic form on the manifold $M_{N}\subset M.$ For this 2-form to be
nondegenerate on $M_{N},$ we assume that Hessian of $\mathcal{L}_{n}$ equals 
$det\left\Vert \frac{\partial ^{2}\mathcal{L}_{n}[u]}{\partial
u_{n+N+1}\partial u_{n+N+1}}\right\Vert \neq 0$ on $M_{N}.$ The latter
proves the theorem.$\rhd $

Let us consider now the prior given dynamical system (7.1) reduced on the
submanifold $M_{N}\subset M.$ To present it as the vector field\ $d/dt$ on $%
M_{N},$ we need preliminary to represent it as a Hamiltonian flow on $M_{N}.$
To do this, let us write the following identities on $M:$ 
\begin{eqnarray}
i_{\frac{d}{dt}}d &<&\mathop{\rm grad}\,\mathcal{L}_{n},du_{n}\rangle =-%
\frac{d}{dn}i_{\frac{d}{dt}}\Omega _{n}^{(2)}[u],  \notag \\
di_{\frac{d}{dt}} &<&\mathop{\rm grad}\,\mathcal{L}_{n},du_{n}\rangle =-%
\frac{d}{dn}(dh_{n}^{(t)}[u]),
\end{eqnarray}%
which are valid for all $n\in \mathbf{Z}.$ Adding the last identities in
(7.7), we come to the following one for all $n\in \mathbf{Z}:$ 
\begin{equation}
\frac{d}{dt}<\mathop{\rm grad}\,\mathcal{L}_{n},du_{n}\rangle =-\frac{d}{dn}%
(i_{\frac{d}{dt}}\Omega _{n}^{(2)}[u]+h_{n}^{(t)}[u]).
\end{equation}%
Having reduced the identity (7.8) upon the manifold $M_{N}\subset M,$ we
obtain the wanted expression for all $u\in M_{N},\quad N\in \mathbf{Z}:$ 
\begin{equation}
i_{\frac{d}{dt}}\Omega _{n}^{(2)}[u]+h_{n}^{(t)}[u]=0.
\end{equation}%
The latter means that the dynamical system (7.1) on the manifold $M_{N}$ is
a Hamiltonian one, with the function $h_{n}^{(t)}[u]$ being a Hamiltonian
one defined explicitly by the second identity in (7.7).

We assume now that the symplectic structure $\Omega _{n}^{(2)}[u]$ on $M_{N}$
be representable as follows: 
\begin{equation}
\Omega _{n}^{(2)}[u]=\sum_{j=0}^{N}\langle dp_{j+n}\wedge du_{j+n}\rangle ,
\end{equation}%
where generalized coordinates $p_{j+n}\in \mathbb{R}^{m},\quad j=\overline{%
0,N},$ are determined from the following discrete jet-expression $\mathcal{L}%
_{n}[u]:=\mathcal{L}(u_{n},u_{n+1},\dots ,u_{n+N+1}),\quad n\in \mathbf{Z},$ 
\begin{equation}
\begin{split}
\alpha _{n}^{(1)}[u]:=\sum_{j=0}^{N}\langle p_{j+n},du_{j+n}\rangle
=\sum_{k=0}^{N}\sum_{j=0}^{k}<\frac{\partial \mathcal{L}_{n+j-k}}{\partial
u_{n+j}},du_{n+j}>  \notag \\
=\sum_{j=0}^{N}\sum_{k=j}^{N}<\frac{\partial \mathcal{L}_{n+j-k}}{\partial
u_{n+j}},du_{n+j}>,
\end{split}%
\end{equation}%
whence we get the final expression: 
\begin{equation}
p_{j+n}:=\sum_{k=j}^{N}\frac{\partial \mathcal{L}_{n+j-k}[u]}{\partial
u_{n+j}},
\end{equation}%
where $j=\overline{0,N},\quad u\in M_{N}\subset M.$

Now we are in a position to reformulate the given dynamical system (7.1) as
that on the reduced manifold $M_{N}\subset M:$ 
\begin{equation}
\frac{du_{n+j}}{dt}=\{h_{n}^{(t)},u_{n+j}\}=\frac{\partial h_{n}^{(t)}}{%
\partial p_{n+j}},\text{ \ \ \ \ \ }\frac{dp_{n+j}}{dt}%
=\{h_{n}^{(t)},p_{n+j}\}=-\frac{\partial h_{n}^{(t)}}{\partial u_{n+j}}
\end{equation}%
for all $n\in \mathbf{Z},j=\overline{0,N}.$ Thereby the problem of embedding
a given discrete dynamical system (7.1) into the vector field\ flow on the
manifold $M_{N}\subset M$ is solved completely with the final result (7.12).

\section{Invariant Lagrangians construction: functional manifold \ case}

\setcounter{equation}{0}

In the case the given nonlinear dynamical system (2.1) being integrable one
of Lax-type, we can proceed effectively to find a commuting infinite
hierarchy of conservation laws can be serving as the invariant Lagrangians,
having been under consideration above.

At first we have to use the important property \cite{Lax} of the
complexified gradient functional $\varphi =\mathop{\rm grad}\,\gamma \in
T^{\ast }(M)\otimes \mathbb{C}$ generated by an arbitrary conservation law $%
\gamma \in D(M),$ i.e. the following Lax-type equation : 
\begin{equation}
d\varphi /dt+{K^{\prime }}^{\ast }\varphi =0
\end{equation}%
where the prime sign denotes the usual Frechet derivative of the local
functional $K:M\rightarrow T(M)$ on the manifold $M,$ the star
\textquotedblright *\textquotedblright\ denotes its conjugation operator
with respect to the nondegenerate standard convolution functional $(\cdot
,\cdot )=\int_{\mathbb{R}}dx\langle \cdot ,\cdot \rangle $ on $T^{\ast
}(M)\times T(M).$ The equation (8.1) admits (what follows from \cite%
{Mit,Nov,Prik}) the special asymptotic solution: 
\begin{equation}
\varphi (x,t;\lambda )\cong (1,a(x,t;\lambda ))^{\tau }exp[\omega
(x,t;\lambda )+\partial ^{-1}\sigma (x,t;\lambda )],
\end{equation}%
where $a(x,t;\lambda )\in \mathbb{C}^{m-1},\quad \sigma (x,t;\lambda )\in 
\mathbb{C},\quad \omega (x,t;\lambda )$ -- some dispersive function. The
sign \textquotedblright $\tau $\textquotedblright\ means here the
transposition one, what is adopted in matrix analysis. For any complex
parameter $\lambda \in \mathbb{C}$ at $|\lambda |\rightarrow \infty $ the
following expansions take place: 
\begin{equation*}
a(x,t;\lambda )\simeq \sum_{j\in \mathbf{Z}_{+}}a_{j}[x,t;u]\lambda
^{-j+s(a)},\text{ }\sigma (x,t;\lambda )\simeq \sum_{j\in \mathbf{Z}%
_{+}}\sigma _{j}[x,t;u]\lambda ^{-j+s(\sigma )}.
\end{equation*}%
Here $s(a)$ and $s(\sigma )\in \mathbf{Z}_{+}$ -- some appropriate
nonnegative integers, the operation $\partial ^{-1}$ means the inverse one
to the differentiation $d/dx,$ that is $d/dx\cdot \partial ^{-1}=1$ for all $%
x\in \mathbb{R}.$

To find the explicit form of the representation (8.2) in the case when the
associated Lax-type representation \cite{Pri} depends parametrically on the
spectral parameter $\lambda (t;z)\in \mathbb{C},$ satisfying the following
non-isospectral condition: 
\begin{equation}
d\lambda (t;z)/dt=g(t;\lambda (t;z)),\qquad \left. \lambda (t;z)\right\vert
_{t=0^{+}}=z\in \mathbb{C},
\end{equation}%
for some meromorphic function $g(t;\cdot ):\mathbb{C}\rightarrow \mathbb{C},$
$t\in \mathbb{R}_{+},$ we must reanalyze more carefully the asymptotic
solutions to the Lax equation (8.1). Namely, we are going to treat more
exactly the case when the solution $\varphi \in T^{\ast }(M)$ to (8.1) is
represented as an appropriate trace-functional of a Lax spectral problem at
the moment $\tau =t\in \mathbb{R}_{+}$ with the spectral parameter $\lambda
(t;\lambda )\in \mathbb{C}$ satisfying the condition (8.3), the evolution of
the given dynamical system (2.1) being considered with respect to the
introduced above parameter $\tau \in \mathbb{R},$ that is 
\begin{equation}
{du}/{d\tau }=K[x,\tau ;u],
\end{equation}%
$\left. u\right\vert _{\tau =0}=\bar{u}\in M$ -- some Cauchy data on $M.$
This means that the functional 
\begin{equation}
\tilde{\varphi}(x,\tau ;\tilde{\lambda}):=reg\mathop{\rm grad}\,Sp\,S(x,\tau
;\tilde{\lambda}),\qquad \tilde{\lambda}=\tilde{\lambda}(\tau ;\lambda
(t;z))\in \mathbb{C},
\end{equation}%
where $S(x,\tau ;\tilde{\lambda})$ is the monodrony matrix corresponding to
a Lax type spectral problem assumed to exist, has to satisfy the
corresponding Lax equation at any point $u\in M$ subject to (8.4): 
\begin{equation}
d\tilde{\varphi}/d\tau +{K^{\prime }}^{\ast }[u]\cdot \tilde{\varphi}=0
\end{equation}%
for all $\tau \in \mathbb{R}_{+}.$ Under the above assumption it is obvious
that the spectral parameter $\tilde{\lambda}=\tilde{\lambda}(\tau ;\lambda
(t;z)),$ where 
\begin{equation}
d\tilde{\lambda}/d\tau =\tilde{g}(\tau ;\tilde{\lambda}),\qquad \left. 
\tilde{\lambda}\right\vert _{\tau =0}=\lambda (t;z)\in \mathbb{C},
\end{equation}%
$\tilde{g}(t;\,\cdot ):\,\mathbb{C}\rightarrow \mathbb{C}$ -- some
meromorphic function being found simply from (8.6) for instance at $u=0,$
the Cauchy data $\lambda (t;z)\in \mathbb{C}$ for all $t\in \mathbb{R}_{+}$
are the corresponding ones to (8.3), the parameter $z\in \mathbb{C}$ being a
spectrum value of the associate Lax type spectral problem at a moment $t\in 
\mathbb{R}_{+}.$

Now we are in a position to formulate the following lemma.

\textsl{Lemma}\textbf{\ 3. }\textit{The Lax equation (8.6) as the parameter }%
$\tau =t\in \mathbb{R}_{+}$\textit{\ admits an asymptotic solution in the
form } 
\begin{equation}
\tilde{\varphi}(x,\tau ;\tilde{\lambda})\cong (1,\tilde{a}(x,\tau ;\tilde{%
\lambda}))^{\tau }exp[\tilde{\omega}(x,\tau ;\tilde{\lambda})+\partial ^{-1}%
\tilde{\sigma}(x,\tau ;\tilde{\lambda})],
\end{equation}%
\textit{where }$\tilde{a}(x,\tau ;\tilde{\lambda})\in \mathbb{C}^{m-1},\quad 
\tilde{\sigma}(x,\tau ;\tilde{\lambda})\in \mathbb{C},$\textit{\ are some
local functionals on }$M,\tilde{\omega}(x,\tau ;\tilde{\lambda})\in \mathbb{C%
}$\textit{\ is\ some dispersion function for all }$x\in \mathbb{R},\tau \in 
\mathbb{R}_{+},$\textit{\ and if \ for }$|\lambda |\rightarrow \infty $%
\textit{\ the property }$|\tilde{\lambda}|\rightarrow \infty $\textit{\ as }$%
\tau =t\in \mathbb{R}_{+}$\textit{\ holds, the following expansions follows: 
} 
\begin{equation}
\tilde{a}(x,\tau ;\tilde{\lambda})\simeq \sum_{j\in \mathbf{Z}_{+}}\tilde{a}%
_{j}[x,\tau ;u]{\tilde{\lambda}}^{-j+s(\tilde{a})},\text{ }\tilde{\sigma}%
(x,\tau ;\tilde{\lambda})\simeq \sum_{j\in \mathbf{Z}_{+}}\tilde{\sigma}%
_{j}[x,\tau ;u]\tilde{\lambda}^{-j+s(\tilde{\sigma})},
\end{equation}%
with $s(\tilde{a})$ and $s(\tilde{\sigma})\in \mathbf{Z}_{+}$ being some
integers.

\textit{Proof.} In virtue of the theory of asymptotic expansions for
arbitrary differential spectral problems, the result (8.8) will hold
provided the representation (8.5) is valid and the spectral parameter $%
\lambda (t;z)\in \mathbb{C}$ is taken subject to (8.7). But the above is the
case because of the Lax-type integrability of the dynamical system (8.4).
Further, due to the mentioned above integrability of (8.4) as well as due to
the well known Stokes property of asymptotic solutions to linear equations
like (8.1), the condition (8.3) holds for some meromorphic function $%
g(t;\cdot ):\mathbb{C}\rightarrow \mathbb{C},$ $t\in \mathbb{R}_{+},$
enjoing the determining property $\frac{d}{dt}\int_{\mathbb{R}}\tilde{\sigma}%
(x,t;\tilde{\lambda}(t;\lambda (t;z)))dx=0$ for all $t\in \mathbb{R}_{+}.$
The latter proves the lemma completely.$\rhd $

As a result of Lemma 3 one ca n formulate the following important theorem.

\textsl{Theorem}\textbf{\ 4.}\textit{\ The Lax integrable parametrically
isospectral dynamical system (8.4) as }$\tau =t\in \mathbb{R}_{+}$\textit{\
admits an infinite hierarchy of conservation laws, in general nonuniform
ones with respect to the variables }$x\in \mathbb{R},\tau \in \mathbb{R}_{+}%
\mathit{,}$\textit{\ which can be represented in an exact form in virtue of
the asymptotic expansion (8.8) and (8.9).}

\textit{Proof. }Indeed, due to the expansion (8.8), we can obtain right away
that the functional 
\begin{equation}
\tilde{\gamma}(t;\lambda (t;z))=\int_{\mathbb{R}}dx\tilde{\sigma}(x,t;\tilde{%
\lambda}(t;\lambda (t;z)))
\end{equation}%
does not depend at $\tau =t\in \mathbb{R}_{+}$ on the parameter $t\in 
\mathbb{R}_{+},$ that is 
\begin{equation}
\left. d\tilde{\gamma}/d\tau \right\vert _{\tau =t\in \mathbb{R}_{+}}=0
\end{equation}%
for all $t\in \mathbb{R}_{+}.$ If we put also the parameter $\tau \in 
\mathbb{R}_{+}$ to tend to $t\in \mathbb{R}_{+},$ due to (8.5) we obtain
that $\left. \tilde{\varphi}(x,\tau ;\tilde{\lambda})\right\vert _{\tau
=t\in \mathbb{R}_{+}}\rightarrow \varphi (x,t;\lambda )$ for all $x\in 
\mathbb{R},t\in \mathbb{R}_{+}$ and $\lambda (t;z)\in \mathbb{C}.$ This
means that a complexified local functional $\varphi (x,t;z)\in T^{\ast
}(M)\otimes \mathbb{C}$ satisfies the equation (8.1) at each point $u\in M.$
As an obvious result, the following identifications hold: 
\begin{equation*}
\left. \tilde{\omega}(x,\tau ;\tilde{\lambda})\right\vert _{\tau =t\in 
\mathbb{R}_{+}}\rightarrow \omega (x,t;z),\text{ \ }\left. \tilde{\sigma}%
(x,\tau ;\tilde{\lambda})\right\vert _{\tau =t\in \mathbb{R}_{+}}\rightarrow
\sigma (x,t;z)
\end{equation*}%
for all $z\in \mathbb{C}.$ Hence, the functional \newline
$\gamma (z):=\left. \tilde{\gamma}(\tau ;\lambda (t;z))\right\vert _{\tau
=t\in \mathbb{R}_{+}}=\int_{\mathbb{R}}dx\sigma (x,t;z)\in D(M)$ doesn't
depend on the evolution parameter $t\in \mathbb{R}_{+}$ and due to equation
(8.1) is a conserved quantity for the nonlinear dynamical system (2.1) under
consideration, i.e. 
\begin{equation}
d\gamma (t;z)/dt=0
\end{equation}%
for all $t\in \mathbb{R}_{+}$ and $z`\in \mathbb{C}.$ Therefore, it makes it
possible to use the equation (8.12) jointly with (8.7) for the asymptotic
expansions (8.9) and (8.3) to be found in exact form. To do this we at first
need to insert the asymptotic expansion (8.8) in the determining equation
(8.6) for the asymptotic expansions (8.9) to be found explicitly at the
moment $\tau =t\in \mathbb{R}_{+}.$ Keeping in mind that at $\tau =t\in 
\mathbb{R}_{+}\quad |\lambda |\rightarrow \infty $ if $|\tilde{\lambda}%
|\rightarrow \infty ,$ and solving step by step the resulting recurrence
relationships for the coefficients in (8.9), we will get the functional $%
\gamma (z):=\left. \tilde{\gamma}(\tau ;\lambda (t;z))\right\vert _{\tau
=t\in \mathbb{R}_{+}},\quad z\in \mathbb{C},$ in the form fitting for the
criteria equation (8.12) could be used. As the second step, we need to use
the differential equation (8.7) for the criteria equation (8.12) to be
satisfied point-wise for all $t\in \mathbb{R}_{+}.$ This means, in
particular, that 
\begin{equation}
\frac{d\gamma (z)}{dt}=\frac{d}{dt}(\left. \sum_{j\in \mathbf{Z}_{+}}\int_{%
\mathbb{R}}dx\tilde{\sigma}_{j}[x,\tau ;u]\tilde{\lambda}^{-j+s(\tilde{\sigma%
})}\right\vert _{\tau =t\in \mathbb{R}_{+}})
\end{equation}%
\begin{equation*}
\left. =\int_{\mathbb{R}}dx\sum_{j\in \mathbf{Z}_{+}}\left[ \frac{d\tilde{%
\sigma}_{j}[x,\tau ;u]}{dt}\tilde{\lambda}^{-j+s(\tilde{\sigma})}+\tilde{%
\sigma}_{j}[x,\tau ;u]\tilde{\lambda}^{-j+s(\tilde{\sigma})-1}(s(\tilde{%
\sigma})-j)\frac{d\tilde{\lambda}}{dt}\right] \right\vert _{\tau =t\in 
\mathbb{R}_{+}}\Rightarrow
\end{equation*}%
\begin{equation*}
\Rightarrow \int_{\mathbb{R}}dx\sum_{j\in \mathbf{Z}_{+}}\left[ \left. (d%
\tilde{\sigma}_{j}/dt)\tilde{\lambda}^{-j+s(\tilde{\sigma}%
)}+\sum_{k>>-\infty }(s(\tilde{\sigma})-k)\tilde{\sigma}_{k}\tilde{g}%
_{j-k-1}(t)\tilde{\lambda}^{-j+s(\tilde{\sigma})}\right\vert _{\tau =t\in 
\mathbb{R}_{+}}+\right.
\end{equation*}%
\begin{equation*}
\left. +\left. \sum_{j\in \mathbf{Z}_{+}}\tilde{\sigma}_{j}[x,t;u]\tilde{%
\lambda}^{-j+s(\tilde{\sigma})-1}(s(\tilde{\sigma})-j)\frac{\partial \tilde{%
\lambda}}{\partial \lambda }g(t;\lambda )\right\vert _{\tau =t\in \mathbb{R}%
_{+}}\right] \equiv 0,
\end{equation*}%
where we have put by definition $\tilde{g}(\tau ;\tilde{\lambda}):\simeq
\sum_{k>>-\infty }\tilde{g}_{k}(\tau ){\tilde{\lambda}}^{-k}$ for $\tau \in 
\mathbb{R}_{+}$ and $|{\tilde{\lambda}}|\rightarrow \infty .$ Since the
spectral parameter $\lambda =\lambda (t;z)$ at the moment $t=0^{+}$
coincides with an arbitrary complex value $z\in \mathbb{C},$ the condition $%
|z|\rightarrow \infty $ together with (8.13) at the moment $t=0^{+}$ gives
rise to the following recurrent relationships: 
\begin{equation}
\sum_{j\in \mathbf{Z}_{+}}\left. \left[ \partial {\tilde{\sigma}}_{j}/dt+{%
\tilde{\sigma}^{\prime }}_{j}\cdot K[t;u]+\sum_{k>>-\infty }(s(\tilde{\sigma}%
)-k){\tilde{\sigma}}_{k}\cdot \tilde{g}_{j-k-1}\right] \tilde{\lambda}^{-j+s(%
\tilde{\sigma})}\right\vert _{\tau =t\in \mathbb{R}_{+}}=
\end{equation}%
\begin{equation*}
=\left. \sum_{j\in \mathbf{Z}_{+}}\tilde{\sigma}_{j}(s(\tilde{\sigma})-j)%
\frac{\partial \tilde{\lambda}}{\partial \lambda }g(t;\lambda )\tilde{\lambda%
}^{-j+s(\tilde{\sigma})-1}\right\vert _{\tau =t\in \mathbb{R}_{+}}\equiv 0%
\text{ }\ (mod\text{ }d/dx)
\end{equation*}%
for all $j\in \mathbf{Z}_{+},x\in \mathbb{R},t\in \mathbb{R}_{+}$ and $u\in
M.$ Having solved the algebraic relationships (8.14) for the prior unknown
function $g(t;\lambda ),$ $t\in \mathbb{R}_{+},$ we will obtain the
generating functional $\gamma (z),$ $z\in \mathbb{C},$ of conservation laws
for (2.1) in exact form. This completes the constructive part of the proof
of the theorem above.

From the practical point of view we need first to get the differential
equation (8.7) in exact, maybe in asymptotic form and find further the
dispersive function $\tilde{\omega}(x,t;\tilde{\lambda})$ and the local
generating functional $\tilde{\sigma}(x,\tau ;\tilde{\lambda})$ defined via
(8.8) and (8.9) for all $x\in \mathbb{R},\tau \in \mathbb{R}_{+}$ and $|%
\tilde{\lambda}|\rightarrow \infty ,$ and next one can find the equation
(8.3) due to$\ $the a$\lg $orithm based on the relationships (8.14). This
together with the possibility of applying the general scheme of the
gradient-holonomic algorithm \cite{PM} gives rise to determining in many
cases the above mentioned Lax-type representation completely in exact form,
that successfully solves the pretty complex direct problem of the
integrability theory of nonlinear dynamical systems on functional manifolds.

Having obtained the generation function $\gamma (z)\in D(M),$ $\ z\in 
\mathbb{C},$ of an infinite hierarchy of conservation laws of the dynamical
system (2.1) on the manifold $M,$ we can build appropriately a general
Lagrangian functional $\mathcal{L}_{N}\in D(M)$ as follows: 
\begin{equation}
\mathcal{L}_{N}=-\gamma _{N+1}+\sum_{j=0}^{N}c_{j}\gamma _{j},
\end{equation}%
where, by definition, $\gamma (z)=\int_{\mathbb{R}}dx\sigma (x,t;z)$ and for 
$|z|\rightarrow \infty $ functionals $\gamma _{j}=\int_{\mathbb{R}}dx\sigma
_{j}[x,t;z],$ $j\in \mathbf{Z}_{+},$ are conservation laws due to expansion
(8.2), with $c_{j}\in \mathbb{R},$ $j=\overline{0,N,}$ being\ some arbitrary
constants and $N\in \mathbf{Z}_{+}$ being an arbitrary but fixed nonnegative
integer. If the differential order of the functional $\gamma _{N+1}\in D(M)$
has the highest one of the orders of functionals $\gamma _{j}\in D(M),j=%
\overline{0,N},$ and additionally, this Lagrangian is not degenerate, that
is $det(Hess\,\gamma _{N+1})\neq 0,$ we can apply in general amost all the
theory developed before, to prove that the critical submanifold $%
M_{N}=\{u\in M:\mathop{\rm grad}\,\mathcal{L}_{N}=0\}$ is a
finite-dimensional invariant manifold inserted into the standard
jet-manifold $J^{(\infty )}(\mathbb{R};\mathbb{R}^{m})$ with the canonical
symplectic structure subject to which our dynamical system is a
finite-dimensional Hamiltonian flow on the invariant submanifold $M_{N}.$

\section{Invariant Lagrangian construction : discrete manifold case}

\setcounter{equation}{0}

Let us consider the discrete Lax integrable dynamical system on the discrete
manifold $M$ without an a priory given Lax-type representation. The problem
arises how to get the corresponding conservation laws via the
gradient-holonomic algorithm \cite{Pri}. To realize this way let us study
solutions to the Lax equation: 
\begin{equation}
d\varphi _{n}/dt+K_{n}^{\prime }[\tau ,u]\cdot \varphi _{n}=0
\end{equation}%
local functionals $\varphi _{n}[u]\in T_{u_{n}}^{\ast }(M)$ at the point $%
u_{n}\in M,n\in \mathbf{Z}.$ In analogy to the approach of Chapter 7 we
assert that equation (9.1) admits a comlexified generating solution $\varphi
_{n}=\varphi _{n}(t;\lambda )\in T_{u_{n}}^{\ast }(M)\otimes \mathbb{C},n\in 
\mathbf{Z},$ $\ $with $z\in \mathbb{C}$ being a complex parameter in the
following form: 
\begin{equation}
\varphi _{n}(t;z)\cong (1,a_{n}(t;z))^{\tau }exp[\omega (t;z)]\biggl(%
\prod_{j=-\infty }^{n}\sigma _{j}(t;z)\biggr),
\end{equation}%
where $\omega (t;z)$ is some dispersive function for $t\in \mathbb{R}%
_{+},\quad a_{n}(t;z)\in \mathbb{C}^{m-1},$\newline
$\sigma _{n}(t;z)\in \mathbb{R}$ are local functionals on $M,$ having the
following asymptotic expansions at $|z|\rightarrow \infty $: 
\begin{equation}
a_{n}(t;z)\simeq \sum_{j\in \mathbf{Z}_{+}}a_{n}[t;u]z^{-j+s(a)},\text{ \ }%
\sigma _{n}^{(j)}(t;z)\simeq \sum_{j\in \mathbf{Z}_{+}}\sigma
_{n}[t;u]z^{-j+s(\sigma )}.
\end{equation}%
To find the explicit form of the asymptotic representation (9.2) we need to
study additionally the asymptotic solutions to the following attached Lax
equation with respect to the new evolution parameter $\tau \in \mathbb{R}%
_{+} $: 
\begin{equation}
d\tilde{\varphi}_{n}/d\tau +K{_{n}^{\prime }}^{\ast }[\tau ,u]\cdot \tilde{%
\varphi}_{n}=0,
\end{equation}%
where $\tilde{\varphi}_{n}\in T_{u_{n}}^{\ast }(M)\otimes \mathbb{C},$ and a
point $u\in M$ evolves subject to the following dynamical system: 
\begin{equation}
du_{n}/d\tau =K_{n}[\tau ;u],\quad
\end{equation}%
for all $n\in \mathbf{Z.}$ \ Having made the assumption above we can assert
based on the general theory of asymptotic solutions to linear equations like
(9.4), that it admits also in general another asymptotic solution in the\
similar form: 
\begin{equation}
\tilde{\varphi}_{n}(\tau ;\tilde{\lambda})\cong (1,\tilde{a}_{n}(\tau ;%
\tilde{\lambda}))^{\tau }exp[\tilde{\omega}(\tau ;\tilde{\lambda}%
)]\prod_{j=-\infty }^{n}\tilde{\sigma}_{j}(\tau ;\tilde{\lambda}),
\end{equation}%
where for all $n\in \mathbf{Z}$ \ and at $\tau \in \mathbb{R}_{+}$ \ the
asymptotic expansions 
\begin{eqnarray}
\tilde{a}_{n}(\tau ;\tilde{\lambda}) &\simeq &\sum_{j\in \mathbf{Z}_{+}}{%
\tilde{a}_{n}}^{(j)}[x,\tau ;u]\tilde{\lambda}^{-j+s(\tilde{a})},  \notag \\
\tilde{\sigma}_{n}(\tau ;\tilde{\lambda}) &\simeq &\sum_{j\in \mathbf{Z}_{+}}%
{\tilde{\sigma}_{n}}^{(j)}[\tau ;u]\tilde{\lambda}^{-j+s(\tilde{\sigma})}
\end{eqnarray}%
hold. The expansions above are valid if $|\tilde{\lambda}|\rightarrow \infty 
$ as $|\lambda (t;z)|\rightarrow \infty ,$ $z\in \mathbb{C}.$ The latter is
the case because of the Lax-integrability of the dynamical system (9.5). The
evolution 
\begin{equation}
d\tilde{\lambda}/d\tau =\tilde{g}(\tau ;\tilde{\lambda}),\qquad \left. 
\tilde{\lambda}\right\vert _{t=0}=\lambda (t;z)\in \mathbb{C},
\end{equation}%
where $\tilde{g}(\tau ;\cdot ):\mathbb{C}\rightarrow \mathbb{C}$ is some
meromorphic mapping for all $\tau \in \mathbb{R}_{+},$ is in general found
making use of the corresponding solution to (9.4) at $u=0.$

Substituting the expansions (9.6) and (9.7) into (9.4), we obtain some
recurrence relationships, giving rise to a possibility the expressions for
local functionals $\tilde{\sigma}_{j}[t;u_{n}],j\in \mathbf{Z}_{+},$ which
can be found exactly. Having this done successfully, we assert that the
functional 
\begin{equation}
\gamma (t;z)=\sum_{n\in \mathbf{Z}}\left. \ln \tilde{\sigma}_{n}(\tau ;%
\tilde{\lambda})\right\vert _{\tau =t\in \mathbb{R}_{+}}\Rightarrow
\sum_{n\in \mathbf{Z}}\ln \sigma _{n}(t;z),
\end{equation}%
where $\tilde{\lambda}=\tilde{\lambda}(\tau ;\lambda ),\,\tau \in \mathbb{R}%
_{+},\,$\ and $\lambda (t;z)\in \mathbb{C},$ is a meromorphic solution to
equation 
\begin{equation}
d\lambda /dt=g(t;\lambda ),\left. \lambda \right\vert _{t=0^{+}}=z\in 
\mathbb{C}
\end{equation}
with still independent meromorphic function $g(t;\cdot )$ for almost all $%
t\in \mathbb{R}_{+}.$ The latter can be found making use of the following
determining condition: the local functional $\left. \tilde{\varphi}_{n}(\tau
;\tilde{\lambda})\right\vert _{\tau =t\in \mathbb{R}_{+}}\rightarrow \varphi
_{n}(t;z)\in T^{\ast }(M)\otimes \mathbb{C}$ for all $t\in \mathbb{R}_{+}$
and $z\in \mathbb{C}.$ Hence, the following equality holds immediately: 
\begin{equation}
\frac{d}{dt}(\sum_{n\in \mathbf{Z}_{+}}\left. \ln \tilde{\sigma}_{n}(\tau ;%
\tilde{\lambda})\right\vert _{\tau =t\in \mathbb{R}_{+}})=
\end{equation}%
\begin{equation*}
\sum_{n\in \mathbf{Z}_{+}}\left. {\tilde{\sigma}}_{n}^{-1}(t;\tilde{\lambda})%
\left[ \frac{\partial \tilde{\sigma}_{n}}{\partial t}+{\tilde{\sigma}}%
_{n}^{\prime }\cdot K_{n}[u]+\frac{\partial \tilde{\sigma}_{n}}{\partial 
\tilde{\lambda}}g(t;\tilde{\lambda})+\right. \right\vert _{\tau =t\in 
\mathbb{R}_{+}}
\end{equation*}%
\begin{equation*}
\left. +\left. \frac{\partial \tilde{\sigma}_{n}}{\partial \tilde{\lambda}}%
\frac{\partial \tilde{\lambda}}{\partial \lambda }g(t;\lambda )\right\vert
_{\tau =t\in \mathbb{R}_{+}}\right] =0
\end{equation*}%
for all $t\in \mathbb{R}_{+}.$ Equating coefficients of (9.11) at all powers
of the spectral parameter $\lambda (t;z)\in \mathbb{C}$ to zero modulus $%
d/dn,$ $n\in \mathbf{Z},$ we will find the recurrent relationships for the
function $g(t;\lambda )$ of (9.8) to be determined successfully. Thereby,
using the equation (9.10) and an expansion $\sigma (t;z)\simeq \sum_{j\in 
\mathbf{Z}_{+}}\sigma _{j}[t;u_{n}]z^{-j+s(\gamma )}$ for $|z|\rightarrow
\infty ,$ where $s(\sigma )\in \mathbf{Z}_{+}$ is some integer number, we
obtain an infinite hierarchy of discrete-wise conservation laws of the
initially given nonlinear dynamical system (2.1) on the manifold $M.$ But
because of the parametric dependence of the conservation laws built above on
the evolution parameter $t\in \mathbb{R}_{+},$ we cannot use right now the
theory developed before to prove the Hamiltonian properties of the
corresponding vector fields on the invariant submanifolds. To do this in an
appropriate way, it is necessary to augment the theory developed before in
some important details.

\section{The reduction procedure on nonlocal Lagrangian submanifolds}

\setcounter{equation}{0}

\noindent \textbf{1. The general algebraic scheme.} Let $\tilde{\mathcal{G}}%
:=C^{\infty }(\mathbb{S}^{1};\mathcal{G})$ be a Lie algebra of loops, taking
values in a matrix Lie algebra $\mathcal{G}$. By means of $\tilde{\mathcal{G}%
}$ one constructs the Lie algebra $\hat{\mathcal{G}}$ of matrix
integral-differential operators \cite{PSA}: 
\begin{equation}
\hat{a}:=\sum_{j\ll \infty }a_{j}\xi ^{j},
\end{equation}%
where the symbol $\xi :=\partial /\partial x$ signs the differentiation with
respect to the independent variable $x\in \mathbb{R}/2\pi \mathbf{Z}\simeq 
\mathbb{S}^{1}$. The usual Lie commutator on $\hat{\mathcal{G}}$ is defined
as: 
\begin{equation}
\lbrack \hat{a},\hat{b}]:=\hat{a}\circ \hat{b}-\hat{b}\circ \hat{a}
\end{equation}%
for all $\hat{a},\hat{b}\in \hat{\mathcal{G}}$, where "$\circ $" is the
product of integral-differential operators, taking the form: 
\begin{equation}
\hat{a}\circ \hat{b}:=\sum_{\alpha \in \mathbf{Z}_{+}}\frac{1}{\alpha !}%
\frac{\partial ^{\alpha }\hat{a}}{\partial \xi ^{\alpha }}\frac{\partial
^{\alpha }\hat{b}}{\partial x^{\alpha }}.
\end{equation}%
On the Lie algebra $\hat{\mathcal{G}}$ there exists the $ad$invariant
nondegenerate symmetric bilinear form: 
\begin{equation}
(\hat{a},\hat{b}):=\int_{0}^{2\pi }Tr\,(\hat{a}\circ \hat{b})\ dx,
\end{equation}%
where $Tr$-operation for all $\hat{a}\in \hat{\mathcal{G}}$ is given by the
expression: 
\begin{equation}
Tr\,\hat{a}:=res_{\xi }\,Sp\,\hat{a}=Sp\,a_{-1},
\end{equation}%
with $Sp$ being the usual matrix trace. With the scalar product (10.4) the
Lie algebra $\hat{\mathcal{G}}$ is transformed into a metrizable one. As a
consequence, the dual to $\hat{\mathcal{G}}$ \ linear space of the matrix
integral-differential operators $\hat{\mathcal{G}}^{\ast }$ is naturally
identified with the Lie algebra $\hat{\mathcal{G}},$ that is $\hat{\mathcal{G%
}}^{\ast }\simeq \hat{\mathcal{G}}$.

The linear subspaces $\hat{\mathcal{G}}_{+}\subset \hat{\mathcal{G}}$ and $%
\hat{\mathcal{G}}_{-}\subset \hat{\mathcal{G}}$ such as 
\begin{eqnarray}
\hat{\mathcal{G}}_{+}\!\!:= &&\!\!\left\{ \hat{a}:=\sum_{j=0}^{n(\hat{a})\ll
\infty }a_{j}\xi ^{j}:\ a_{j}\in \tilde{\mathcal{G}},j=\overline{0,n(\hat{a})%
}\right\} ,  \notag \\
\hat{\mathcal{G}}_{-}\!\!:= &&\!\!\left\{ \hat{b}:=\sum_{j=0}^{\infty }\xi
^{-(j+1)}b_{j}:\ b_{j}\in \tilde{\mathcal{G}},j\in \mathbf{Z}_{+}\right\} ,
\end{eqnarray}%
are Lie subalgebras in $\hat{\mathcal{G}}$ and $\hat{\mathcal{G}}=\hat{%
\mathcal{G}}_{+}\oplus \hat{\mathcal{G}}_{-}$. Because of the splitting of $%
\hat{\mathcal{G}}$ into the direct sum of its Lie subalgebras one can
construct the so called Lie-Poisson structure \cite{Ad, Bl, Oe, PM} on $\hat{%
\mathcal{G}}^{\ast },$ using a special linear endomorphism $\mathcal{R}$ of $%
\hat{\mathcal{G}}$: 
\begin{equation}
\mathcal{R}:=(P_{+}-P_{-})/2,\quad P_{\pm }\hat{\mathcal{G}}:=\hat{\mathcal{G%
}}_{\pm },\ \ P_{\pm }\hat{\mathcal{G}}_{\mp }=0.
\end{equation}

For any smooth by Frechet functionals $\gamma ,\mu \in \mathcal{D}(\hat{%
\mathcal{G}}^{\ast })$ the Lie-Poisson bracket on $\hat{\mathcal{G}}^{\ast }$
is given by the expression: 
\begin{equation}
\left\{ \gamma ,\mu \right\} _{\mathcal{R}}(\hat{l})=\left( \hat{l},[\nabla
\gamma (\hat{l}),\nabla \mu (\hat{l})]_{\mathcal{R}}\right) ,
\end{equation}%
where $\hat{l}\in \hat{\mathcal{G}}^{\ast }$ and for all $\hat{a},\ \hat{b}%
\in \hat{\mathcal{G}}$ the $\mathcal{R}$commutator in (10.8) has the form 
\cite{Oe, PM}: \ \ \ \ \ \ \ \ \ \ 
\begin{equation}
\lbrack \hat{a},\hat{b}]_{\mathcal{R}}:=[\mathcal{R}\hat{a},\hat{b}]+[\hat{a}%
,\mathcal{R}\hat{b}],
\end{equation}%
subject to which the linear space $\hat{\mathcal{G}}$ becomes a Lie algebra
too. The gradient $\nabla \gamma (\hat{l})\in \hat{\mathcal{G}}$ of a
functional $\gamma \in \mathcal{D}(\hat{\mathcal{G}}^{\ast })$ at a point $%
\hat{l}\in \hat{\mathcal{G}}^{\ast }$ with respect to the scalar product
(10.4) is defined as 
\begin{equation}
\delta \gamma (\hat{l}):=\left( \nabla \gamma (\hat{l}),\delta \hat{l}%
\right) ,
\end{equation}%
where the linear space isomorphism $\hat{\mathcal{G}}\simeq \hat{\mathcal{G}}%
^{\ast }$ is taken into account.

The Lie-Poisson bracket (10.8) generates Hamiltonian dynamical systems on $%
\hat{\mathcal{G}}^{\ast }$ related with Casimir invariants $\gamma \in I(%
\mathcal{G}^{\ast })$, satisfying the condition: 
\begin{equation}
\lbrack \nabla \gamma (\hat{l}),\hat{l}]=0,
\end{equation}%
as the corresponding Hamiltonian functions. Due to the expressions (10.8)
and (10.11)the mentioned above Hamiltonian system takes the form: 
\begin{equation}
d\hat{l}/dt:=[\mathcal{R}\nabla \gamma (\hat{l}),\hat{l}]=[\nabla \gamma
_{+}(\hat{l}),\hat{l}],
\end{equation}%
being equivalent to the usual commutator Lax type representation \cite{PM,
La}. The relationship (10.12) is a compatibility condition for the linear
integral-differential equations: 
\begin{eqnarray}
\hat{l}f\!\! &=&\!\!\lambda f,  \notag \\
df/dt\!\! &=&\!\!\nabla \gamma _{+}(\hat{l})f,
\end{eqnarray}%
where $\lambda \in \mathbb{C}$ is a spectral parameter and a vector-function 
$f\in W(\mathbb{S}^{1};\mathbf{H})$ is an element of some matrix
representation for the Lie algebra $\hat{\mathcal{G}}$ in some functional
Banach space $\mathbf{H}$.

Algebraic properties of the equation (10.12) together with (10.14) and the
associated dynamical system on the space of adjoint functions $f^{\ast }\in
W^{\ast }(\mathbb{S}^{1};\mathbf{H})$: 
\begin{equation}
df^{\ast }/dt=-(\nabla \gamma (\hat{l}))_{+}^{\ast }f^{\ast },
\end{equation}%
where $f^{\ast }\in W^{\ast }$ is a solution to the adjoint spectral
problem: 
\begin{equation}
\hat{l}^{\ast }f^{\ast }=\nu f^{\ast },
\end{equation}%
being considered as some coupled evolution equations on the space $\hat{%
\mathcal{G}}^{\ast }\oplus W\oplus W^{\ast }$ is an object of our further
investigation.

\textbf{2. The tensor product of Poisson structures and its Backlund
transformation.} To compactify the description below we will use the
following designation of the gradient vector 
\begin{equation*}
\nabla \gamma (\tilde{l},\tilde{f},\tilde{f}^{\ast }):=(\delta \gamma
/\delta \tilde{l},\,\delta \gamma /\delta \tilde{f},\,\delta \gamma /\delta 
\tilde{f}^{\ast })^{T}
\end{equation*}%
for any smooth functional $\gamma \in \mathcal{D}(\hat{\mathcal{G}}^{\ast
}\oplus W\oplus W^{\ast })$. On the spaces $\hat{\mathcal{G}}^{\ast }$ and $%
W\oplus W^{\ast }$ there exist canonical Poisson structures \cite{Bl, OS, PM}
\begin{equation}
\delta \gamma /\delta \tilde{l}:\overset{\tilde{\theta}}{\rightarrow }%
[(\delta \gamma /\delta \tilde{l})_{+},\tilde{l}]-[\delta \gamma /\delta 
\tilde{l},\tilde{l}]_{+}
\end{equation}%
at a point $\tilde{l}\in \hat{\mathcal{G}}^{\ast }$ and 
\begin{equation}
(\delta \gamma /\delta \tilde{f},\,\delta \gamma /\delta \tilde{f}^{\ast
})^{T}:\overset{\tilde{J}}{\rightarrow }(\delta \gamma /\delta \tilde{f}%
^{\ast },\,-\delta \gamma /\delta \tilde{f})^{T}
\end{equation}%
at a point $(\tilde{f},\tilde{f}^{\ast })\in W\oplus W^{\ast }$
correspondingly. It should be noted that the Poisson structure (10.17) is
transformed into (10.12) for any Casimir functional $\gamma \in I(\hat{%
\mathcal{G}}^{\ast })$ . Thus, on the extended space $\hat{\mathcal{G}}%
^{\ast }\oplus W\oplus W^{\ast }$ one can obtain a Poisson structure as the
tensor product $\tilde{\Theta}:=\tilde{\theta}\otimes \tilde{J}$ of the
structures (10.17) and (10.18).

Let us consider the following Backlund transformation \cite{PM, OS, SP}: 
\begin{equation}
(\hat{l},f,f^{\ast }):\overset{B}{\rightarrow }(\tilde{l}(\hat{l},f,f^{\ast
}),\tilde{f}=f,\tilde{f}^{\ast }=f^{\ast }),
\end{equation}%
generating on $\hat{\mathcal{G}}^{\ast }\oplus W\oplus W^{\ast }$ a Poisson
structure $\Theta $ with respect to variables $(\hat{l},f,f^{\ast })$ of the
coupled evolution equations (10.12), (10.14), (10.15).

The main condition for the mapping (10.19) to be defined is the coincidence
of the dynamical system 
\begin{equation}
(d\hat{l}/dt,\,df/dt,\,df^{\ast }/dt)^{T}:=-\Theta \nabla \gamma (\hat{l}%
,f,f^{\ast })
\end{equation}%
with (10.12), (10.14), (10.15) in the case of $\gamma \in I(\hat{\mathcal{G}}%
^{\ast }),$ i.e. if this functional is taken to be not dependent of
variables $(f,f^{\ast })\in W\oplus W^{\ast }$. To satisfy that condition,
one has to find a variation of any smooth Casimir functional $\gamma \in I(%
\hat{\mathcal{G}}^{\ast })$ as $\delta \tilde{l}=0,$ considered as a
functional on $\hat{\mathcal{G}}^{\ast }\oplus W\oplus W^{\ast },$ taking
into account flows (10.14), (10.15) and the Backlund transformation (10.19): 
\begin{eqnarray}
\left. \delta \gamma (\tilde{l},\tilde{f},\tilde{f}^{\ast })\right\vert
_{\delta \tilde{l}=0}=(<\delta \gamma /\delta \tilde{f},\delta \tilde{f}%
>)+(<\delta \gamma /\delta \tilde{f}^{\ast },\delta \tilde{f}^{\ast }>)\!\!
&=&\!\!  \notag \\
\left. (<-d\tilde{f}^{\ast }/dt,\delta \tilde{f}>)+(<d\tilde{f}/dt,\delta 
\tilde{f}^{\ast }>)\right\vert _{\tilde{f}=f,\,\tilde{f}^{\ast }=f^{\ast
}}\!\! &=&\!\!  \notag \\
(<(\delta \gamma /\delta \hat{l})_{+}^{\ast }f^{\ast },\delta f>)+(<(\delta
\gamma /\delta \hat{l})_{+}f,\delta f^{\ast }>)\!\! &=&\!\!  \notag \\
(<f^{\ast },(\delta \gamma /\delta \hat{l})_{+}\delta f>)+(<(\delta \gamma
/\delta \hat{l})_{+}f,\delta f^{\ast }>)\!\! &=&\!\!  \notag \\
(\delta \gamma /\delta \hat{l},\delta f\xi ^{-1}\otimes f^{\ast })+(\delta
\gamma /\delta \hat{l},f\xi ^{-1}\otimes \delta f^{\ast })\!\! &=&\!\! 
\notag \\
(\delta \gamma /\delta \hat{l},\delta (f\xi ^{-1}\otimes f^{\ast
})):=(\delta \gamma /\delta \hat{l},\delta \hat{l}).\!\! &\ &\!\!
\end{eqnarray}%
As a result of the expression (10.21) one obtains the relationships: 
\begin{equation}
\left. \delta \hat{l}\right\vert _{\delta \tilde{l}=0}=\delta (f\xi
^{-1}\otimes f^{\ast }),
\end{equation}%
or having assumed the linear dependence of $\hat{l}$ and $\tilde{l}\in \hat{%
\mathcal{G}}^{\ast }$ one gets right away that 
\begin{equation}
\hat{l}=\tilde{l}+f\xi ^{-1}\otimes f^{\ast }.
\end{equation}%
Thus, the Backlund transformation (10.19) can be now written as 
\begin{equation}
(\hat{l},f,f^{\ast }):\overset{B}{\rightarrow }(\tilde{l}=\hat{l}-f\xi
^{-1}\otimes f^{\ast },\text{ \ }\tilde{f}=f,\text{ \ }\tilde{f}^{\ast }=%
\tilde{f}^{\ast }).
\end{equation}%
The expression (10.24) generalizes the result, obtained in the papers \cite%
{PM, SP} for the Lie algebra $\hat{\mathcal{G}}$ of integral-differential
operators with scalar coefficients. The existence of the Backlund
transformation (10.19) makes it possible to formulate the following theorem.

\textbf{Theorem 1.} \textit{A dynamical system on $\hat{\mathcal{G}}^{\ast
}\oplus W\oplus W^{\ast }$, being Hamiltonian with respect to the canonical
Poisson structure $\tilde{\Theta}:T^{\ast }(\hat{\mathcal{G}}^{\ast }\oplus
W\oplus W^{\ast })\rightarrow T(\hat{\mathcal{G}}^{\ast }\oplus W\oplus
W^{\ast })$, and generated by the evolution equations: 
\begin{equation}
d\tilde{l}/dt=[\nabla \gamma _{+}(\tilde{l}),\tilde{l}]-[\nabla \gamma (%
\tilde{l}),\tilde{l}]_{+},\quad d\tilde{f}/dt=\delta \gamma /\delta \tilde{f}%
^{\ast },\quad d\tilde{f}^{\ast }/dt=-\delta \gamma /\delta \tilde{f},
\end{equation}%
with $\gamma \in I(\mathcal{G}^{\ast })$ being the Casimir functional at $%
\hat{l}\in \hat{\mathcal{G}}^{\ast }$ connected with $\tilde{l}\in \hat{%
\mathcal{G}}^{\ast }$ by (10.23), is equivalent to the system (10.12),
(10.14) and (10.15) via the constructed above Backlund transformation
(10.24).}

By means of simple calculations via the formula (see for egz. \cite{PM, Bl}) 
\begin{equation*}
\tilde{\Theta}=B^{^{\prime }}\Theta B^{^{\prime }\ast },
\end{equation*}%
where $B^{^{\prime }}:T(\hat{\mathcal{G}}^{\ast }\oplus W\oplus W^{\ast
})\rightarrow T(\hat{\mathcal{G}}^{\ast }\oplus W\oplus W^{\ast })$ is the
Frechet derivative of (24), one brings about the following form of the
Poisson structure $\Theta $ on $\hat{\mathcal{G}}^{\ast }\oplus W\oplus
W^{\ast }\ni (\hat{l},f,f^{\ast })$: 
\begin{equation*}
\nabla \gamma (\hat{l},f,f^{\ast }):\overset{\Theta }{\rightarrow }%
\begin{pmatrix}
\left[ \hat{l},(\delta \gamma /\delta \hat{l})_{+}\right] -\left[ \hat{l}%
,\delta \gamma /\delta \hat{l}\right] _{+}- \\ 
-f\xi ^{-1}\otimes \delta \gamma /\delta f+\delta \gamma /\delta f^{\ast
}\xi ^{-1}\otimes f^{\ast }\ \delta \gamma /\delta f^{\ast }-(\delta \gamma
/\delta \hat{l})_{+}f\ -\delta \gamma /\delta f+(\delta \gamma /\delta \hat{l%
})_{+}^{\ast }f%
\end{pmatrix}%
\end{equation*}%
that makes it possible to formulate the next theorem.

\textbf{Theorem 2.} \textit{The dynamical system (10.20), being Hamiltonian
with respect to the Poisson structure $\Theta $ in the form (10.26) and a
functional $\gamma \in I(\hat{\mathcal{G}}^{\ast })$, gives the inherited
Hamiltonian representation for the coupled evolution equations (10.12),
(10.14), (10.15).}

By means of the expression (10.23) one can construct Hamiltonian evolution
equations, describing commutative flows on the extended space $\hat{\mathcal{%
G}}^{\ast }\oplus W\oplus W^{\ast }$ at a fixed element $\tilde{l}\in \hat{%
\mathcal{G}}^{\ast }.$ Due to (10.24) every equation of such a type is
equivalent to the system 
\begin{equation}
\left\{ 
\begin{array}{l}
d\hat{l}/d\tau _{n}=[\hat{l}_{+}^{n},\hat{l}], \\ 
df/d\tau _{n}=\hat{l}_{+}^{n}f, \\ 
df^{\ast }/d\tau _{n}=-(\hat{l}^{\ast })_{+}^{n}f^{\ast },%
\end{array}%
\right.
\end{equation}%
generated by involutive with respect to the Poisson bracket (10.17) Casimir
invariants $\gamma _{n}\in I(\hat{\mathcal{G}}^{\ast }),$ $n\in \mathbf{N},$
taking here the standard form: 
\begin{equation*}
\gamma _{n}=1/(n+1)(\hat{l}^{n},\hat{l})
\end{equation*}%
at $\hat{l}\in \hat{\mathcal{G}}^{\ast }$.

The compatibility conditions of the Hamiltonian systems (10.25) for
different $n\in \mathbb{Z}_{+}$ can be used for obtaining Lax integrable
equations on usual spaces of smooth $2\pi $-periodic multivariable functions
that will be done in the next section.

\textbf{3. The Lax type integrable Davey-Stewartson equation and its triple
linear representation.} Choose the element $\tilde{l}\in \hat{\mathcal{G}}%
^{\ast }$ in an exact form such as 
\begin{equation*}
\tilde{l}=\left( 
\begin{array}{cc}
1 & 0 \\ 
0 & -1%
\end{array}%
\right) \xi -\left( 
\begin{array}{cc}
0 & u \\ 
\bar{u} & 0%
\end{array}%
\right) ,
\end{equation*}%
where $u,\bar{u}\in C^{\infty }(\mathbb{S}^{1};\mathbb{C})$ and $\mathcal{G}%
=gl(2;\mathbb{C})$. Then 
\begin{equation}
\hat{l}=\tilde{l}+\left( 
\begin{array}{cc}
f_{1}\xi ^{-1}f_{1}^{\ast } & \ f_{1}\xi ^{-1}f_{2}^{\ast }+u \\ 
\bar{u}+f_{2}\xi ^{-1}f_{1}^{\ast } & f_{2}\xi ^{-1}f_{2}^{\ast }%
\end{array}%
\right) ,
\end{equation}%
where $f=(f_{1},f_{2})^{T}$ and $f^{\ast }=(f_{1}^{\ast },f_{2}^{\ast
})^{T}, $ "$^{-}$" can sign the complex conjugation. Below we will study the
evolutions (10.25) of vector-functions $(f,f^{\ast })\in W(\mathbb{S}^{1};%
\mathbb{C}^{2})\oplus W^{\ast }(\mathbb{S}^{1};\mathbb{C}^{2})$ with respect
to the variables $y=\tau _{1}$ and $t=\tau _{2}$ at the point (10.26). They
can be obtained from the second and third equations in (10.25), having put $%
n=1$ and $n=2$, as well as from the first one. The latter is the
compatibility condition of the spectral problem 
\begin{equation}
\hat{l}\Phi =\lambda \Phi ,
\end{equation}%
where $\Phi =(\Phi _{1},\Phi _{2})^{T}\in W(\mathbb{S}^{1};\mathbb{C}^{2})$, 
$\lambda \in \mathbb{C}$ is some parameter, with the following linear
equations: 
\begin{eqnarray}
d\Phi /dy\!\! &=&\!\!\hat{l}_{+}\Phi , \\
d\Phi /dt\!\! &=&\!\!\hat{l}_{+}^{2}\Phi ,
\end{eqnarray}%
arising from (10.26) at $n=1$ and $n=2$ correspondingly. The compatibility
of equations (10.28) and (10.29) leads to the relationships: 
\begin{eqnarray}
&&\partial u/\partial y=-2f_{1}f_{2}^{\ast },\quad \partial \bar{u}/\partial
y=-2f_{1}^{\ast }f_{2}, \\
&&\partial f_{1}/\partial y=\partial f_{1}/\partial x-uf_{2},\quad \partial
f_{1}^{\ast }/\partial y=\partial f_{1}^{\ast }/\partial x-\bar{u}%
f_{2}^{\ast },  \notag \\
&&\partial f_{2}/\partial y=-\partial f_{2}/\partial x+\bar{u}f_{1},\quad
\partial f_{2}^{\ast }/\partial y=-\partial f_{2}^{\ast }/\partial
x+uf_{1}^{\ast }.  \notag
\end{eqnarray}%
Analogously, replacing $t\in \mathbb{C}$ by $it\in i\mathbb{R},\ i^{2}=-1$,
one gets from (10.29) and (10.30): 
\begin{eqnarray}
&&\!\!\!\!\!\!\!\!\!\!\!\!\!\!\!\!\!\!\!\!du/dt=i(\partial ^{2}u/\partial
x\partial y+2u(f_{1}f_{1}^{\ast }+f_{2}f_{2}^{\ast })),\quad d\bar{u}%
/dt=-i(\partial ^{2}\bar{u}/\partial x\partial y+2\bar{u}(f_{1}f_{1}^{\ast
}+f_{2}f_{2}^{\ast })),  \notag \\
&&\!\!\!\!\!\!\!\partial (f_{1}f_{1}^{\ast })/\partial y-\partial
(f_{1}f_{1}^{\ast })/\partial x=1/2\partial (u\bar{u})/\partial y=-(\partial
(f_{2}f_{2}^{\ast })/\partial x+\partial (f_{2}f_{2}^{\ast })/\partial y), 
\notag \\
&&\quad \quad df_{1}/dt=i(\partial ^{2}f_{1}/\partial
x^{2}+(2f_{1}f_{1}^{\ast }-u\bar{u})f_{1}-\partial u/\partial xf_{2}), 
\notag \\
&&\quad \quad df_{1}^{\ast }/dt=-i(\partial ^{2}f_{1}^{\ast }/\partial
x^{2}+(2f_{1}f_{1}^{\ast }-u\bar{u})f_{1}^{\ast }-\partial \bar{u}/\partial
xf_{2}^{\ast }), \\
&&\quad \quad df_{2}/dt=i(\partial ^{2}f_{2}/\partial
x^{2}-(2f_{2}f_{2}^{\ast }+u\bar{u})f_{2}-\partial \bar{u}/\partial xf_{1}),
\notag \\
&&\quad \quad df_{2}^{\ast }/dt=-i(\partial ^{2}f_{2}^{\ast }/\partial
x^{2}-(2f_{2}f_{2}^{\ast }+u\bar{u})f_{2}^{\ast }-\partial u/\partial
xf_{1}^{\ast }).  \notag
\end{eqnarray}%
The relationships (10.31), (10.32) take the well known form of the
Davey-Stewartson equation \cite{Bl, KSS} at $\bar{u}\in $C$^{\infty }(%
\mathbb{S}^{1};\mathbb{C})$ being a complex conjugated to $u\in $C$^{\infty
}(\mathbb{S}^{1};\mathbb{C})$. The compatibility for every pair of equations
(10.28), (10.29) and (10.30), which can be rewritten as the first order
linear ordinary differential ones in such a way: 
\begin{eqnarray}
&&d\Phi /dx=\left( 
\begin{array}{ccc}
\lambda & u & -f_{1} \\ 
\bar{u} & -\lambda & f_{2} \\ 
f_{1}^{\ast } & f_{2}^{\ast } & 0%
\end{array}%
\right) \Phi , \\
&&d\Phi /dy=\left( 
\begin{array}{ccc}
\lambda & 0 & -f_{1} \\ 
0 & \lambda & -f_{2} \\ 
f_{1}^{\ast } & f_{2}^{\ast } & 0%
\end{array}%
\right) \Phi , \\
&&d\Phi /dt=i\left( 
\begin{array}{ccc}
\lambda ^{2}+f_{1}f_{1}^{\ast } & 1/2\partial u/\partial y & -\lambda
f_{1}-\partial f_{1}/\partial y \\ 
-1/2\partial \bar{u}/\partial y & \lambda ^{2}-f_{2}f_{2}^{\ast } & -\lambda
f_{2}-\partial f_{1}/\partial y \\ 
\lambda f_{1}^{\ast }+\partial f_{1}^{\ast }/\partial y & \lambda
f_{2}^{\ast }+\partial f_{2}^{\ast }/\partial y & 0%
\end{array}%
\right) \Phi ,
\end{eqnarray}%
where $\Phi =(\Phi _{1},\Phi _{2},\Phi _{3})^{T}\in W(\mathbb{S}^{1};\mathbb{%
C}^{3}),$ provide its Lax type integrability. Thus, the following theorem
holds.

\textbf{Theorem 3.} \textit{The Davey-Stewartson equation (10.32), (10.33)
possesses the Lax representation as the compatibility condition for
equations (10.34) under the additional natural constraint (10.27).}

In fact, one has found above a triple linearization for a (2+1)-dimensional
dynamical system, that is a new important ingredient of the Lie algebraic
approach to Lax type integrable flows, based on the Backlund type
transformation (\textit{10.}23) developed in this work. It is clear that the
similar construction of a triple linearization like (\textit{10.}4) can be
done for many other both old and new (2+1)-dimensional dynamical systems, on
what we plant to stop in detail in another work under preparation.

\section{Conclusion}

\setcounter{equation}{0}

The developed above theory of parametrically Lax-type integrable dynamical
systems concedes to widen to a great extent the class of exactly treated
nonlinear models in many fields of science. It is to be noted here the
following important mathematical fact being got in the paper: almost every
nonlinear dynamical system admits a parametrically isospectral Lax type
representation but a given dynamical system is the Lax-type integrable if an
evolution of the spectrum parameter doesn't depend on a point $u\in M$ at
all Cauchy data. This result has allowed us to develop a very effective
direct criterion for the following problem: whether a given nonlinear
dynamical system on the functional manifold $M$ is parametrically Lax-type
integrable or not. Having the problem above solved, we have suggested the
reduction procedure for the associated nonlinear dynamical systems to be
descended on the invariant submanifold $M_{N}\subset M$ built before
inheriting the canonical Hamiltonian structure and the Liouville complete
integrability. Thereby, the powerful techniques of perturbation theory can
be successfully used for dynamical systems under consideration, as well as
the relationships between the full Hamiltonian theory and various
Hamiltonian truncations could be now got understandable more deeply.

The imbedding problem for infinite-dimensional dynamical systems with
additional structures such as invariants and symmetries is as old as the
Newton-Lagrange mechanics, having been treated by many researches, using
both analytical and algebraic methods. The powerful differential-geometric
tools used here were created mainly in works by E.~Cartan at the beginning
of the twentieth century. The great impact in the development of imbedding
methods was done in last time, especially owing to theory of isospectral
deformations for some linear structures built on the special vector bundles
over the spase $M$ as the base of a given nonlinear dynamical system. Among
them there are such structures as the moment map $l:M\rightarrow \mathcal{G}%
^{\ast }$ into the adjoint space to the Lie algebra $\mathcal{G}$ of
symmetries, acting on the symplectic phase space $M$ equivariantly \cite%
{Pri,Gill}, the connection of the Cartan-Eresman structures appearing via
the Wahlquist-Estabrook approach \cite{Wah}, and many others.

For the last years the general structure of Lagrangian and Hamiltonian
formalisms was studied thoroughly using both geometrical and algebraical
methods \cite{Kup,Kupe}. The special attention was paid to the theory of
differential-difference dynamical systems on the infunite-dimensional
manifolds \cite{Kupe,Deift}. Some number of articles was devoted to the
theory of pure discrete dynamical systems \cite{Mos,Bae,Ves,Levi}, as well
treating the interesting examples \cite{Levi} appeared to be important for
applications.

In future work we intend to treat further imbedding problems for
infinite-dimensional both continuous and discrete dynamical systems basing
on the differential-geometric Cartan's theory of differential ideals in
Grassmann algebras over jet-manifolds, intimately connected with the problem
under regard. As it is well known, there existed by now only two regular
enough algorithmic approaches \cite{PSA, PM, SP} to constructing integrable
multi-dimensional (mainly 2+1) dynamical systems on infinite-dimensional
functional spaces. Our approach, devised in this work, is substantially
based on the results previously done in \cite{PM, SP}, explains completely
the computational properties of multi-dimensional flows before delivered in
works \cite{Bl, KSS}. As the key points of our approach there used the
canonical Hamiltonian structures naturally existing on the extended phase
space and the related with them Backlund transformation which saves Casimir
invariants of a chosen matrix integral-differential Lie algebra. The latter
gives rise to some additional Hamiltonian properties of considered extended
evolution flows before studied in \cite{Bl, PM} making use of the standard
inverse scattering transform \cite{La, Bl, PM} and the formal symmetry
reduction for the KP-ierarchy \cite{OS, SP} of commuting operator flows.

As one can convince ourselves analyzing the structure of the Backlund type
transformation (10.24), that it strongly depends on the type of an $ad$%
-invariant scalar product chosen on an operator Lie algebra $\hat{\mathcal{G}%
}$ and its Lie algebra decomposition like (10.6). Since there exist in
general other possibilities of choosing such decompositions and $ad$%
-invariant scalar products on $\hat{\mathcal{G}}$, they give rise naturally
to another resulting types of the corresponding Backlund transformations,
which can be a subject of another special investigation. Let us here only
mention the choice of a scalar product related with the operator Lie algebra 
$\hat{\mathcal{G}}$ centrally extended by means of the standard
Maurer-Cartan two-cocycle \cite{PSA, Oe, PM}, bringing about new types of
multi-dimensional integrable flows.

The last aspect of the Backlund approach to constructing Lax type integrable
flows and their partial solutions which is worth of mention is related with
Darboux-Backlund type transformations \cite{Bl, MS} and their new
generalization recently developed in \cite{Ni, SP}. They give rise to very
effective procedures of constructing multi-dimensional integrable flows on
functional spaces with arbitrary number of independent variables
simultaniously delivering a wide class of their exact analytical solutions,
depending on many constant parameters, which can appear to be useful for
diverse applications in applied sciences.

All mentioned above Backlund type transformations aspects can be studied as
special investigations, giving rise to new directions in the theory of
multi-dimensional evolution flows and their integrability.

\section{\protect\bigskip Acknowledgements}

The authors thank Profs. D.L. Blackmore (NJIT, Newark, NJ,USA), M.O.
Perestiuk (Kyiv National University, Ukraina) for useful comments on the
results of the work.


\bigskip

\begin{thebibliography}{10}
\bibitem[1]{Grif} P. A. Griffiths. Exterior Differential Systems and the
Calculus of Variations. New York. Birkhause, 1982, 480p.

\bibitem[2]{Fil} Prykarpatsky A. K., B. M. Fil. Category of topological
jet-manifolds and certain applications in the theory of nonlinear
infinite-dimensional dynamical systems. // Ukr. Math. J., 1993, vol. 44, No
2 p. 1136-1147.

\bibitem[3]{Lax} P. D. Lax. Periodic solutions of the Korteweg-de Vries
equation. // Comm. Pure and Appl. Mathem., 1968, vol. 21, No 2, p. 467-490.

\bibitem[4]{Bry} R. L. Bryant. On notions of Equivalence of Variational
Problems with One Independent Variable. // Contemp. Math., 1987, vol. 63, p.
65-76.

\bibitem[5]{St} S. Sternberg. Some preliminary remarks on the formal
variational calculus of Gel'fand and Dikii. // Lect. Notes in Mathem., 1980,
vol. 150, p. 399-407.

\bibitem[6]{Stern} S. Sternberg. Lectures on differential geometry. Moscow,
Mir, 1961, p. 412.

\bibitem[7]{Bogo} O. Bogoyavlensky, S. Novikov. On connection of Hamiltonian
formalisms of stationar and nonstationar problems. // Function. Anal.

\bibitem[8]{Pri} Prykarpatsky A. K., I. V. Mykytiuk. Algebraic aspects of
integrable dynamical systems on manifolds. Kiev, Nauk. dumka, 1991, 380p.

\bibitem[9]{Gelf} I. M. Gelfand, L. A. Dikiy Integrable nonlinear equations
and Liouville theorem. // Function. Anal. Appl., 1979, vol. 13, No 1, p.
8-20.

\bibitem[10]{Gill} V.Gillemin, S.Sternberg. The moment map and Collective
motion. // Ann. Phys., 1980, vol. 127, No 2, p. 220-253

\bibitem[11]{Wah} Wahlgnist H. D., Estabrook F. B. Prolongation structures
of nonlinear evolution equations. // J. Math. Phys., 1975, vol. 16, No 1, p.
1-7, 1976, vol. 17, No 7, p. 1293-1297.

\bibitem[12]{Kup} Kupershmidt B. A. Geometry of jet-bundles and the
structure of Lagrangian and Hammmiltonian formalisms. // Lect. Notes Math.,
1980, vol. 775, p. 162-218.

\bibitem[13]{Kupe} Kupershmidt B. A. Discrete Lax equations and
differential-difference calculus. // Asterisque, 1985, vol. 123, p. 5-212.

\bibitem[14]{Duft} Deift P. Li L.-C., Tomei C. Loop groups, discrete
versions of some classical integrable systems and rank-2 extensions. //
Memoirs of the AMS, 1992, vol. 100, No 479, p. 1-101

\bibitem[15]{Mos} Moser J., Veselov A. P. Discrete versions of some
classical integrable systems and factorization of matrix polynomials.
Zurich, ETH, Preprint, 1989, 76p.

\bibitem[16]{Bae} Baez J. C., Gillam J. W. An algebraic approach to discrete
mechanics. // Lett. Math. Phys., 1994, vol. 31, No 3, p. 205-212.

\bibitem[17]{Ves} Veselov A. P. What is an itegrable mapping? in: What is
itegrability? Ney York, Springer-Verlag, 1991, p. 251-272.

\bibitem[18]{Levi} Levi D., Winternitz P. Continuous symmetries of discrete
equations. // Physics Lett. A., 1991, vol. 152, No 7, p. 335-338.

\bibitem[19]{Kan} Kaneko K. Symplectic cellular automate. // Phys. Lett. A.,
1988, vol. 129, No 1, p. 9-16.

\bibitem[20]{Mars} Albraham R., Marsden J. Foundation of mechanics. London,
The Blujamin Publish. Co., 1978, p. 806.

\bibitem[21]{Sou} Souriau J. M. Structure des systemes dynamique. Paris,
Dunod, 1970.

\bibitem[22]{Adl} Adler M. On a trace functional for formal
pseudo-differential operators and the symplectic structures of the
Korteweg-de Vries equations. // Invent. Mathem., 1979, vol. 50, No 2, p.
219-248.

\bibitem[23]{Pry} Prykarpatsky A.K. and others. Algebraic structure of the
gradient-golonomic algorithm for Lax integrable nonlinear dynamical systems.
// J. Math. Phys., vol. 35, No 4, p. 1763-1777, vol. 35, No 8, p. 6115-6126.

\bibitem[24]{Oew} Oewel W. Dirac constraints in field theory: Lifts of
Hamiltonian systems to the cotangent bundle. // J. Math. Phys., 1988, vol.
29, No 1, p. 210-219.

\bibitem[25]{Oe} Oewel W. R-structures, Yang-Baxter equations and related
involution theorems. // J. Math. Phys., 1989, vol. 30, No 5, p. 1140-1149.

\bibitem[26]{Fok} Fokas A. S., Gelfand I. M. Bi-Hamiltonian Structures and
Integrability. in: Important developments in Soliton Theory.
Springer-Verlag, 1992.

\bibitem[27]{Olv} Olver P. J. Canonical forms and integrability of
bi-Hamiltonian systems. // Phys. Lett. A., 1990, vol. 148, No 3, p. 177-187

\bibitem[28]{Mag} Magri F. A simple model of the integrable Hamiltonian
equation. // J. Math. Phys., 1978, vol. 19, No 3, p. 1156-1162.

\bibitem[29]{Fer} Fernandes R. L. Completely integrable bi-Hamiltonian
systems. // J. of Dynam. and Diff. Equ., 1994, vol. 6, No 1, p. 53-69.

\bibitem[30]{Mit} Mitropolsky Yu. A., Bogoliubov N. N., Prykarpatsky A. K.,
Samoilenko V.H. Integrable dynamical systems. Spectral and
differential-geometric aspects. Kiev, Naukova Dumka, 1987.

\bibitem[31]{Nov} S.P.Novikov (editor). The theory of solitons. Moscov, Mir,
1980.

\bibitem[32]{Prik} Mitropolsky Yu. O., Prykarpatsky A. K., Fil B. M. Some
aspects of a gradient-holonomic algorithm in the theory of integrability of
nonlinear dynamical systems and computer algebra problems // Ukrainian Math.
J., 1991, vol. 43, No 1, p. 63-74.

\bibitem[33]{Kuy} Kuybida V. V., Prytula M. M., Prykarpatsky A. K. The study
of properties of the parametric isospectral integrability of nonlinear
dynamical systems on functional manifolds and their finitedimensional
approximations. Preprint of Inst. for Appl. Problems of Mech. and Mathem. of
Ukr. Academy of Sci., Lviv, 1991, No 10, 41p.

\bibitem[34]{Wein} Marsden J., Weinstein A. Reduction of symplectic
manifolds with symmetries. // Rep. Math. Phys., 1974, vol. 5, No 2, p.
121-130.

\bibitem[35]{Str} Oewel W., Strampp W. Constrained KP-hierarchy and
bi-Hamiltonian structures. // Comm. Math. Phys., 1993, vol. 157, No 1, p.
51-81.

\bibitem[36]{PSA} Prykarpatsky A.K., Samoilenko V.Hr., Andrushkiw R.I.,
Mitropolsky Yu.O., Prytula M.M.\} Algebraic structure of the
gradient-holonomic algorithm for Lax integrable nonlinear systems. I // J.
Math. Phys. -- 1999. -- 35, {2116} 4. -- {421}. 1763-1777.

\bibitem[37]{Ad} Adler M. On a trace functional for formal
pseudo-differential operators and the symplectic structures of a Korteweg-de
Vries equation // Invent. Math. -- 1979. -- 50, {2116} 2. -- P. 219-248.

\bibitem[38]{Bl} Blaszak M. Multi-Hamiltonian theory of dynamical systems.
-- Verlag-Berlin-Heidelberg: Springer, 1998. -- 345 p.

\bibitem[39]{Oe} Oevel W. \ R-structures, Yang-Baxter equations and related
involution theorems // J. Math. Phys. -- 1989. -- 30, {2116} 5. -- P.
1140-1149.

\bibitem[40]{PM} Prykarpatsky A.K., Mykytiuk I.V. Algebraic integrability of
nonlinear dynamical systems on manifolds: classical and quantum aspects. --
Dordrecht-Boston-London: Kluwer Academic Publishers, 1998. -- 553 p.

\bibitem[41]{La} Lax P.D. Periodic solutions of the KdV equation // Commun.
Pure and Appl. Math. -- 1975. -- \textbf{28}. -- {421}. 141-188.

\bibitem[42]{OS} Oevel W., Strampp W. Constrained KP hierarchy and
bi-Hamiltonian structures // Commun. Math. Phys. -- 1993. -- 157. -- P.
51-81.

\bibitem[43]{KSS} Konopelchenko B., Sidorenko Yu., Strampp W.
(1+1)-dimensional integrable systems as symmetry constraints of
(2+1)-dimensional systems // Phys. Lett. A. -- 1991. -- 157. -- {421}. 17-21.

\bibitem[44]{MS} Matveev V.B., Salle M.I. Darboux-Backlund transformations
and applications. -- New York: Springer, 1993.

\bibitem[45]{Ni} Nimmo J.C.C. In: Nonlinear evolution equations and
dynamical systems (NEEDS'94) / Ed. by Makhankov V.G., Bishop A.R. and Holm
D.D. -- World Scient. Publ., 1994.

\bibitem[46]{SP} Samoilenko A.M., Prykarpatsky Ya.A. Algebraic-analytic
aaspects of integrable nonlinear dynamical systems and their perturbations.
-- Kyiv: Intitute of Mathematics at NAS of Ukraine, 2002
\end{thebibliography}
\end{document}